\documentclass[times,sort&compress,3p]{elsarticle}
\journal{Journal of Multivariate Analysis}
\usepackage[labelfont=bf]{caption}
\usepackage{amsmath,amsfonts,amssymb,amsthm,booktabs,color,epsfig,graphicx,hyperref,url}
\usepackage{pstricks}
\usepackage{float}
\theoremstyle{plain}
\newtheorem{theorem}{Theorem}

\newtheorem{lemma}{Lemma}

\newtheorem{assumption}{Assumption}
\theoremstyle{definition}

\begin{document}
\begin{frontmatter}
\title{Nonparametric drift estimation from diffusions with correlated brownian motions}
\author[1]{Fabienne COMTE}
\author[2]{Nicolas MARIE\corref{mycorrespondingauthor}}
\address[1]{Universit\'e Paris Cit\'e, CNRS, MAP5 UMR 8145, F-75006 Paris, France}
\address[2]{Laboratoire Modal'X, Universit\'e Paris Nanterre, Nanterre, France}
\cortext[mycorrespondingauthor]{Corresponding author. Email address: \url{n.marie@parisnanterre.fr}}
\begin{abstract}
In the present paper, we consider that $N$ diffusion processes $X^1,\dots,X^N$ are observed on $[0,T]$, where $T$ is fixed and $N$ grows to infinity. Contrary to most of the recent works, we no longer assume that the processes are independent. The dependency is modeled through correlations between the Brownian motions driving the diffusion processes. A nonparametric estimator of the drift function, which does not use the knowledge of the correlation matrix, is proposed and studied. Its integrated mean squared risk is bounded and an adaptive procedure is proposed. Few theoretical tools to handle this kind of dependency are available, and this makes our results new. Numerical experiments show that the procedure works in practice.
\end{abstract}
\end{frontmatter}
%

% Section : Introduction.

%
\section{Introduction}\label{introduction_section}
We start by describing our model. Consider the diffusion process $X = (X_t)_{t\in [0,T]}$, defined by
\begin{equation}\label{main_equation}
X_t = x_0 +\int_{0}^{t}b(X_s)ds +
\int_{0}^{t}\sigma(X_s)dB_s
\textrm{ $;$ }
t\in [0,T],
\end{equation}
where $x_0\in\mathbb R$, $B = (B_t)_{t\in [0,T]}$ is a Brownian motion, $b :\mathbb R\rightarrow\mathbb R$ is a Lipschitz continuous function, and $\sigma :\mathbb R\rightarrow\mathbb R$ is a bounded Lipschitz continuous function. Now, let $(B^1,\dots,B^N)$ be a $N$-dimensional Gaussian process ($N\in\mathbb N/\{0\}$) whose components are copies of $B$ such that
\begin{equation}\label{dependence_condition_B}
\mathbb E(B_{s}^{i}B_{t}^{k}) =
R_{i,k}(s\wedge t)
\textrm{ $;$ }
\forall i,k\in\{1,\dots,N\}\textrm{, }
\forall s,t\in [0,T],
\end{equation}
where $R = (R_{i,k})_{i,k}$ is a correlation matrix. Note that, thanks to the (stochastic) integration by parts formula, the dependence condition (\ref{dependence_condition_B}) on $B^1,\dots,B^N$ implies that, for every $i,k\in\{1,\dots,N\}$, $d\langle B^i,B^k\rangle_t = R_{i,k}dt$, with $R_{i,i} = 1$. Finally, consider $X^i :=\mathcal I(x_0,B^i)$ for every $i\in\{1,\dots,N\}$, where $\mathcal I(.)$ is the It\^o map associated to Equation (\ref{main_equation}). In the present paper, we consider that these  $N$ diffusion processes $X^1,\dots,X^N$ are observed on $[0,T]$, where $T$ is fixed and $N$ grows to infinity, and our aim is to estimate nonparametrically the drift function $b(.)$.
\\
\\
In the case of independent Brownian motions, that is $R =\mathbf I_N$ (the $N\times N$ identity matrix), projection least squares estimators have been studied in Comte and Genon-Catalot \cite{CGC20} for continuous time observations, in Denis {\it et al.} \cite{DDM21} for discrete time (with small step) observations with a classification purpose in the parametric setting, and in Denis {\it et al.} \cite{DDM20} in the nonparametric context, for instance. Marie and Rosier \cite{MR22} propose a kernel based Nadaraya-Watson estimator of the drift function $b$, with bandwidth selection relying on the Penalized Comparison to Overfitting (PCO) criterion recently introduced by Lacour {\it et al.} \cite{LMR17}. Still in the case $R =\mathbf I_N$, Comte and Marie \cite{CM21} investigate the properties of the projection least squares estimator of the drift when $B$ is a fractional Brownian motion.\\
Dependency is often encountered in recent works in the context of stochastic systems of $N$ interacting particles, with recent nonparametric drift estimators proposals in Della Maestra and Hoffmann \cite{DMH22}, Belomestny {\it et al.} \cite{BPP21} or Comte and Genon-Catalot \cite{CGC22}. These kinds of models are related to physics. We rather have in mind economic or financial models. For instance, in Duellmann {\it et al.} \cite{DKK10}, the authors consider a portfolio of $N$ homogenous firms such that the asset value $X_{t}^{i}$ at time $t$ of the $i$-th firm is modeled by Merton's model (see \cite{MERTON74}) $dX_{t}^{i} =\mu X_{t}^{i}dt +\sigma X_{t}^{i}dB_{t}^{i}$ with $X_{0}^{i} = X_0$, which corresponds to (\ref{main_equation}) with $b(x) =\mu x$ and $\sigma(x) =\sigma x$. Intending to capture the dependency between the firms, they also assume that $dB_{t}^{i} =\sqrt\rho dW_t +\sqrt{1 -\rho}dW_{t}^{i}$, where $W$ is a {\it common systematic risk factor}, $W^i$ is a {\it firm-specific risk factor} and $\rho\in [0,1]$. This corresponds to a particular matrix $R$, precisely $R_{i,k} =\rho$ for $i\neq k$ (and $R_{i,i} = 1$), so that one single parameter $\rho$ represents the so-called {\it asset correlation}. This model has been considered in e.g. Bush {\it et al.} \cite{BHH11}, for the more mathematical purpose of studying the limit of the empirical distribution of the $X_{t}^{i}$'s (see also references therein). Our extension from specific geometric Brownian motion to general nonparametric diffusion (\ref{main_equation}), and from one single correlation parameter to a general matrix representation, is therefore standard in both respects. This context has nevertheless never been considered before up to our knowledge. Let us emphasize that our aim is not to estimate $R$, but to exhibit conditions on it such that $b(.)$ can be estimated with performance near of the independent setting.
\\
\\
In our framework, $T$ is fixed, and $N$ is large. Our results are nonasymptotic, but the idea is that $N$ grows to infinity. We fix a subset $I$ of $\mathbb R$ and build a collection of projection least squares estimators of $b_I = b\mathbf 1_I$ where $I$ is compact or not. The estimators are defined by their coefficients on an orthonormal basis of $\mathbb L^2(I)$, $\varphi_1,\dots,\varphi_m$, resulting from a standard least squares computation. Precisely, we consider the estimator of the drift function $b$ minimizing the objective function $\gamma_N(\tau)$
\begin{equation}\label{objective_function}
\tau\longmapsto \gamma_N(\tau):=\frac{1}{NT}\sum_{i = 1}^{N}\left(
\int_{0}^{T}\tau(X_{s}^{i})^2ds - 2\int_{0}^{T}\tau(X_{s}^{i})dX_{s}^{i}\right)
\end{equation}
on the $m$-dimensional function space $\mathcal S_m = {\rm span}\{\varphi_1,\dots,\varphi_m\}$. The first part of $\gamma_N(\tau)$ involves a quantity
\begin{displaymath}
\|\tau\|_{N}^{2} :=
\frac{1}{NT}\sum_{i = 1}^{N}\int_{0}^{T}\tau(X_{s}^{i})^2ds,
\end{displaymath}
which is considered as {\it the squared empirical norm} of the function $\tau$. These estimators are the same as in Comte and Genon-Catalot \cite{CGC20}, but their study is made significantly more difficult by the dependency context. We do not have at our disposal any coupling method nor any transformation leading to a simpler system; in particular, applying $R^{-1/2}$ to the system does not bring any simplification because of a "widespreading" of the components of the process. Tropp's deviation inequalities used in the independent context (see Tropp \cite{TROPP12}, Matrix Chernov Inequality, Theorem 1.1 and Matrix Bernstein Inequality, Theorem 1.4), which allow to consider the empirical norm and its expectation (an integral norm, thus) as equivalent with high probability, no longer apply. Martingale properties still are useful, and we turn to Azuma's matrix deviation inequality (see Tropp \cite{TROPP12}, Theorem 7.1), which however requires to set  sparsity conditions on $R$ (see Assumption \ref{sparsity_condition}). This equivalence property between empirical and weighted integral norms  is the key of the rigorous study of the risk of the drift estimator, and the correlation matrix is therefore at the heart of the proofs.
\\
\\
The plan of the paper is the following. A first parametric example motivates the model and the way of estimating a drift parameter in Section \ref{motivation_section}. The general nonparametric drift estimator is defined in Section \ref{projection_LS_estimator_section} and a risk bound on a fixed projection space is proved. Adaptive estimation is studied in Section \ref{model_selection_section} and the whole procedure is illustrated through simulations in Section \ref{section_numerical_experiments}. Lastly, proofs are gathered in Section \ref{proofs_section}.
%

% Section : A preliminary simple financial motivation in the parametric framework.

%
\section{Preliminary motivation and example in the parametric framework}\label{motivation_section}
This preliminary section deals  with the geometric model described in the introduction, in the parametric framework, in order to motivate our investigations. Similarly to Duellmann {\it et al.} \cite{DKK10}, consider $N$ risky assets of same nature and of prices processes $X^1,\dots,X^N$ observed on the time interval $[0,T]$. Since these assets are of same nature, to model their prices by $dX_{t}^{i} = \mu X_{t}^{i}dt +\sigma X_{t}^{i}dB_{t}^{i}$ with $\mu\in\mathbb R$ and $\sigma > 0$ not depending on $i\in\{1,\dots,N\}$ is realistic, but it is also very realistic to consider that $B^1,\dots,B^N$ may be dependent, through the correlation matrix described above. Let us compute the quadratic risk of the least squares estimator $\widehat\theta_N$ of $\theta =\mu -\sigma^2/2$ in this special case. Since we can write that $X_{t}^{i} = x_0\exp(Y_{t}^{i})$ with $Y_{t}^{i} =\theta t +\sigma B_{t}^{i}$ for every $i\in\{1,\dots,N\}$ and $t\in [0,T]$, we set
\begin{displaymath}
\widehat\theta_N =
\frac{1}{NT}\sum_{i = 1}^{N}Y_{T}^{i} =
\theta +\frac{\sigma}{NT}\sum_{i = 1}^{N}B_{T}^{i}.
\end{displaymath}
Then,
\begin{eqnarray*}
 \mathbb E(|\widehat\theta_N -\theta|^2) & = &
 \frac{\sigma^2}{N^2T^2}\left[
 \sum_{i = 1}^{N}\mathbb E(|B_{T}^{i}|^2) +
 \sum_{i\neq k}\mathbb E(B_{T}^{i}B_{T}^{k})
 \right]\\
 & = &
 \frac{\sigma^2}{NT} +
 \frac{\sigma^2}{N^2T}\sum_{i\neq k}R_{i,k}
 =\frac{\sigma^2}{NT}\left(1 +\frac{1}{N}\sum_{i\neq k}R_{i,k}\right).
\end{eqnarray*}
This means that the rate of convergence of $\widehat\theta_N$ is of order
\begin{displaymath}
\mathbf V :=\frac{1}{N} +\frac{1}{N^2}\sum_{i\neq k}R_{i,k}.
\end{displaymath}
We note that if $R_{i,k} =\rho$ for all $i\neq k$, then the estimator is not consistent. This would be the same if all the coefficients of $R$ were positive and only bounded by a constant $\rho > 0$. However, if we set a sparsity condition by saying that $R$ is block-diagonal with blocks of size (less than) $k_0$, and if we assume that all nonzero coefficients are equal to (or bounded by) $\rho$, then $\sum_{i\neq k}R_{i,k}\leqslant k_0\rho N$. So, $k_0\rho$ is the loss in risk due to dependency, while the rate remains $O(1/N)$. Referring to the firms model of Duellmann {\it et al.} \cite{DKK10} and Bush {\it et al.}\cite{BHH11}, this means for instance that for a large $N$, dependent firms have to be grouped as several independent sets aggregated in the global model.
\\
\\
Another way to model the dependency with few parameters is to assume that
\begin{displaymath}
dB_{t}^{i} =\sqrt{\tt a}dW_{t}^{i} +\sqrt{1 -\texttt a}dW_{t}^{i + 1},
\end{displaymath}
where $W^1,\dots,W^{N + 1}$ are independent Brownian motions and $\texttt a\in [0,1]$. This is a way of saying that each firm is correlated to the following one in the list. In that case,
\begin{displaymath}
R_{i,i + 1} = R_{i,i - 1} =\sqrt{{\tt a}(1 - {\tt a})},
\quad R_{i,i} = 1,\quad
R_{i,k} = 0\textrm{ for }|k - i| > 1;
\end{displaymath}
and then
\begin{displaymath}
{\bf V} =\frac{1}{N}
\left(1 + 2\left(1 -\frac{1}{N}\right)\sqrt{{\tt a}(1 - {\tt a})}\right)
\end{displaymath}
has order $O(1/N)$. Note that this matrix is sparse in the sense of Assumption \ref{sparsity_condition} below.
\\
\\
Our purpose is to show that, at least for some special dependence schemes on $B^1,\dots,B^N$, the variance term of the projection (nonparametric) least squares estimator of $b(.)$ introduced in the following section is at most of order
\begin{displaymath}
\frac{1}{N}
\left(1 +\frac{1}{N}\sum_{i\neq k}|R_{i,k}|\right).
\end{displaymath}
It is noteworthy that the estimator $\widehat\theta_N$ is the maximum likelihood estimator (MLE) when $X^1,\dots,X^N$ are independent, while the MLE in our dependent setting would involve -- and thus require the complete knowledge of -- the matrix $R$ (more specifically, its inverse). In the present strategy, the knowledge of $R$ is not required, which is interesting and may justify a loss of efficiency.
%

% Section : A projection least squares estimator of the drift function.

%
\section{A projection least squares estimator of the drift function}\label{projection_LS_estimator_section}
%

% Subsection : The objective function.

%
\subsection{The objective function}
Set $N_T := [NT] + 1$ and let $f_T$ be the density function defined by
\begin{displaymath}
f_T(x) :=\frac{1}{T}\int_{0}^{T}p_s(x_0,x)ds
\textrm{ $;$ }
\forall x\in\mathbb R,
\end{displaymath}
where $p_{s}(x_0,.)$ is the density with respect to Lebesgue's measure of the probability distribution of $X_s$ for every $s\in (0,T]$. Let us consider the projection space
$\mathcal S_m := {\rm span}\{\varphi_1,\dots,\varphi_m\}$, where $\varphi_1,\dots,\varphi_{N_T}$ are continuous functions from $I$ into $\mathbb R$ such that $(\varphi_1,\dots,\varphi_{N_T})$ is an orthonormal family in $\mathbb L^2(I,dx)$, and $I\subset\mathbb R$ is a non-empty interval. We recall now that the objective function $\tau\in\mathcal S_m\mapsto\gamma_N(\tau)$ is defined by (\ref{objective_function}), where $m\in\{1,\dots,N_T\}$. We choose a contrast which is the same as in the independent case. Note that for the nonparametric estimation of the drift function from $N$ observed paths, even in the independent case, least squares and maximum likelihood strategies do not match. Indeed, the likelihood would involve weights $\sigma(X_{s}^{i})^{-2}$ inside all integrals. In the dependent case, there would also be the matrix $R^{-1}$ to take into account. Even if both $\sigma(.)$ and $R$ can be considered as known, it is interesting not to need them to compute the drift estimator. In particular, the step to discrete time high frequency data is then much simpler. Since the strategy works in the independent case, we can hope that if the correlations are not too strong, then the strategy remains relevant.
\\
\\
{\bf Remark.} For any $\tau\in\mathcal S_m$,
\begin{eqnarray*}
 \mathbb E(\gamma_N(\tau)) & = &
 \frac{1}{T}\int_{0}^{T}\mathbb E(|\tau(X_s) - b(X_s)|^2)ds -
 \frac{1}{T}\int_{0}^{T}\mathbb E(b(X_s)^2)ds\\
 & = &
 \int_{-\infty}^{\infty}(\tau(x) - b(x))^2f_T(x)dx -
 \int_{-\infty}^{\infty}b(x)^2f_T(x)dx.
\end{eqnarray*}
Then, the more $\tau$ is {\it close} to $b$, the smaller  $\mathbb E(\gamma_N(\tau))$. For this reason, the estimator of $b$ minimizing $\gamma_N(.)$ is studied in this paper.
%

% Subsection : The projection least squares estimator and some related matrices.

%
\subsection{The projection least squares estimator and some related matrices}
In this section, $m$ is a fixed integer in $\{1, \dots, N_T\}$. We consider the estimator
\begin{equation}\label{estimator}
\widehat b_m :=\arg\min_{\tau\in\mathcal S_m}\gamma_N(\tau)
\end{equation}
of $b$, if it exists and is unique. Since $\mathcal S_m = {\rm span}\{\varphi_1,\dots,\varphi_m\}$, there exist $m$ square integrable random variables $\widehat\theta_1,\dots,\widehat\theta_m$ such that
\begin{displaymath}
\widehat b_m =
\sum_{j = 1}^{m}\widehat\theta_j\varphi_j.
\end{displaymath}
Then,
\begin{displaymath}
\nabla\gamma_N(\widehat b_m) =
\left(\frac{1}{NT}\sum_{i = 1}^{N}\left(
2\sum_{\ell = 1}^{m}\widehat\theta_\ell\int_{0}^{T}
\varphi_j(X_{s}^{i})\varphi_{\ell}(X_{s}^{i})ds
- 2\int_{0}^{T}\varphi_j(X_{s}^{i})dX_{s}^{i}\right)\right)_{j\in\{1,\dots,m\}}.
\end{displaymath}
Let 
\begin{displaymath}
\widehat{\bf\Psi}_m :=
\left(\frac{1}{NT}\sum_{i = 1}^{N}\int_{0}^{T}
\varphi_j(X_{s}^{i})\varphi_{\ell}(X_{s}^{i})ds\right)_{j,\ell\in\{1,\dots,m\}}
\end{displaymath}
and
\begin{displaymath}
\widehat{\bf X}_m :=\left(\frac{1}{NT}
\sum_{i = 1}^{N}\int_{0}^{T}
\varphi_j(X_{s}^{i})dX_{s}^{i}\right)_{j\in\{1,\dots,m\}}.
\end{displaymath}
Therefore, by (\ref{estimator}) and if $\widehat{\bf\Psi}_m$ is invertible, necessarily
\begin{displaymath}
\widehat{\bf\Theta}_m :=
(\widehat\theta_1,\dots,\widehat\theta_m)^* =
\widehat{\bf\Psi}_{m}^{-1}\widehat{\bf X}_m,
\end{displaymath}
where $\mathbf M^*$ denotes the transpose of the matrix $\mathbf M$.
\\
\\
{\bf Remarks:}
\begin{enumerate}
 \item We can write $\widehat{\bf\Psi}_m = (\langle\varphi_j,\varphi_{\ell}\rangle_N)_{j,\ell}$, where
 \begin{displaymath}
 \langle\varphi,\psi\rangle_N :=
 \frac{1}{NT}\sum_{i = 1}^{N}\int_{0}^{T}
 \varphi(X_{s}^{i})\psi(X_{s}^{i})ds
 \end{displaymath}
 for every measurable functions $\varphi$ and $\psi$ from $\mathbb R$ into itself.
 \item The following useful decomposition holds: $\widehat{\bf X}_m = (\langle b,\varphi_j\rangle_N)_{j}^{*} +\widehat{\bf E}_m$, where
 \begin{displaymath}
 \widehat{\bf E}_m :=
 \left(\frac{1}{NT}\sum_{i = 1}^{N}\int_{0}^{T}
 \sigma(X_{s}^{i})\varphi_j(X_{s}^{i})dB_{s}^{i}\right)_{j\in\{1,\dots,m\}}^{*}.
 \end{displaymath}
\end{enumerate}
Let us introduce the two following deterministic matrices related to the previous random ones:
\begin{itemize}
 \item ${\bf\Psi}_m :=\mathbb E(\widehat{\bf\Psi}_m) = (\langle\varphi_j,\varphi_{\ell}\rangle_{f_T})_{j,\ell}$, where $\langle .,.\rangle_{f_T}$ is the  scalar product in $\mathbb L^2(I,f_T(x)dx)$.
 \item ${\bf\Psi}_{m,\sigma} := NT\mathbb E(\widehat{\bf E}_m\widehat{\bf E}_{m}^{*})$.
\end{itemize}
Note that under the following assumption, Comte and Genon-Catalot established in \cite{CGC20} (see Lemma 1) that $\mathbf\Psi_m$ is invertible.
%

% Assumption : Assumption on the basis.

%
\begin{assumption}\label{assumption_basis}
The $\varphi_j$'s satisfy the three following conditions:
\begin{enumerate}
 \item $(\varphi_1,\dots,\varphi_m)$ is an orthonormal family of $\mathbb L^2(I,dx)$.
 \item The $\varphi_j$'s are bounded, continuously derivable, and have bounded derivatives.
 \item There exist $x_1,\dots,x_m\in I$ such that $\det[(\varphi_j(x_{\ell}))_{j,\ell}]\not= 0$.
\end{enumerate}
\end{assumption}
\noindent
Let us conclude this section with the following suitable bound on the trace of $\mathbf\Psi_{m}^{-1/2}\mathbf\Psi_{m,\sigma}\mathbf\Psi_{m}^{-1/2}$. To that aim, we define the following quantity associated with the basis:
\begin{displaymath}
L(m) :=
1\vee\left(
\sup_{x\in I}\sum_{j = 1}^{m}\varphi_j(x)^2\right).
\end{displaymath}
%

% Lemma : Preliminary bound on the variance term.

%
\begin{lemma}\label{preliminary_lemma_variance}
Under Assumption \ref{assumption_basis}, for $\sigma$ belonging to $\mathbb L^2(\mathbb R,f_T(x)dx)$ but possibly unbounded,
\begin{equation}\label{trace1}
{\rm trace}(\mathbf\Psi_{m}^{-1/2}\mathbf\Psi_{m,\sigma}\mathbf\Psi_{m}^{-1/2})
\leqslant
\mathfrak c_{\ref{preliminary_lemma_variance}}
L(m)\|\mathbf\Psi_{m}^{-1}\|_{\rm op}
\left(1 +\frac{1}{N}\sum_{i\neq k}|R_{i,k}|\right)
\end{equation}
with
\begin{displaymath}
\mathfrak c_{\ref{preliminary_lemma_variance}} =
\int_{-\infty}^{\infty}\sigma(x)^2f_T(x)dx.
\end{displaymath}
If in addition $\sigma$ is bounded, then
\begin{equation}\label{trace2}
{\rm trace}(\mathbf\Psi_{m}^{-1/2}
\mathbf\Psi_{m,\sigma}\mathbf\Psi_{m}^{-1/2})
\leqslant m\|\sigma\|_{\infty}^{2}
\left(1 +\frac{1}{N}\sum_{i\neq k}|R_{i,k}|\right).
\end{equation}
\end{lemma}
\noindent
The two previous bounds on the trace compete. In some contexts, they can have the same order $m$. This occurs if both following conditions hold (this setting is referred to as "compact setting" below):
\begin{itemize}
 \item $I$ is a compact set and $L(m) = m$, as in the case of a trigonometric basis.
 \item $f_T$ is lower bounded on $I$ by $f_0 > 0$. Indeed, then $\|\mathbf\Psi_{m}^{-1}\|_{\rm op}\leqslant 1/f_0$ (see \cite{CGC20}).
\end{itemize}
However, the bound (\ref{trace1}) is not relevant for a non compactly supported basis. For instance, for the Hermite basis described below, $L(m)$ is of order $\sqrt m$, but $\|\mathbf\Psi_{m}^{-1}\|_{\rm op}$ is increasing with $m$ and can be checked to be numerically very large. So, for the Hermite basis, the second bound (\ref{trace2}) must be preferred if $\sigma$ is bounded, and is used in the sequel.\\
Finally, note that $N^{-1}\sum_{i,k}|R_{i,k}|$ can be replaced by $\|R\|_{\rm op}$ in the bounds (\ref{trace1}) and (\ref{trace2}), and it holds that
\begin{displaymath}
\frac{1}{N}
\left|\sum_{i,k = 1}^{N}R_{i,k}\right|
\leqslant\|R\|_{\rm op}.
\end{displaymath}
Therefore, if the coefficients of the matrix $R$ are nonnegative, as in the examples of Section \ref{motivation_section} or in the simulation Section \ref{section_numerical_experiments}, then $N^{-1}\sum_{i,k}|R_{i,k}|$ is better than $\|R\|_{\rm op}$.
%

% Section : Risk bound on the projection least squares estimator.

%
\subsection{Risk bound on the projection least squares estimator}
This section deals, for a given model $m$, with a risk bound on the truncated estimator
\begin{displaymath}
\widetilde b_m :=\widehat b_m\mathbf 1_{\widehat\Lambda_m},
\end{displaymath}
where
\begin{displaymath}
\widehat\Lambda_m :=
\left\{L(m)(\|\widehat{\bf\Psi}_{m}^{-1}\|_{\rm op}\vee 1)
\leqslant\mathfrak c_T(p)\frac{NT}{\log(NT)}
\right\}
\end{displaymath}
with
\begin{displaymath}
\mathfrak c_T(p) =
\frac{1}{256T(1 + p/2)}
\textrm{, }p\geqslant 12.
\end{displaymath}
On the event $\widehat\Lambda_m$, $\widehat{\bf\Psi}_m$ is invertible because
\begin{displaymath}
\inf\{{\rm sp}(\widehat{\bf\Psi}_m)\}
\geqslant\frac{L(m)}{\mathfrak c_T(p)}\cdot\frac{\log(NT)}{NT} \, >0,
\end{displaymath}
and then $\widetilde b_m$ is well-defined. Consider
\begin{displaymath}
\mathcal I_N :=
\{i\in\{2,\dots,N\} :\exists k\in\{1,\dots,i - 1\}\textrm{ such that }R_{i,k}\neq 0\}.
\end{displaymath}
In the sequel, $m$ fulfills the following assumption, related to the stability condition introduced in Cohen {\it et al.} \cite{CDL13} and also used in Comte and Genon-Catalot \cite{CGC20}. Due to dependency, it has to be reinforced by undesirable squares.
%

% Assumption : Condition on m.

%
\begin{assumption}\label{condition_m}
$\displaystyle{[L(m)(\|\mathbf\Psi_{m}^{-1}\|_{\rm op}\vee 1)]^2
\leqslant\frac{\mathfrak c_T(p)}{2}\cdot\frac{NT}{\log(NT)}}$.
\end{assumption}
\noindent
Note that in the so-called {\it compact setting} defined above, this condition reduces to
\begin{displaymath}
m\lesssim\sqrt{\frac{NT}{\log(NT)}},
\end{displaymath}
which is similar to the constraint obtained in Baraud \cite{BARAUD02} (see the condition $N_n\leqslant K^{-1}\sqrt{n/\log(n)^3}$ in his Theorem 1.1). However, this last condition can be improved in the independent case.
\\
\\
Moreover, a sparsity condition has to be set on $|\mathcal I_N|$, and this is again in order to handle dependency.
%

% Assumption : Sparsity condition.

%
\begin{assumption}\label{sparsity_condition}
There exists a deterministic constant $\mathfrak c_{\ref{sparsity_condition}} > 0$, not depending on $m$ and $N$, such that
\begin{displaymath}
|\mathcal I_N|\leqslant
\mathfrak c_{\ref{sparsity_condition}}N^{\frac{1}{2} -\frac{6}{p}}.
\end{displaymath}
\end{assumption}
\noindent
{\bf Remarks:}
\begin{enumerate}
 \item The dependence condition (\ref{dependence_condition_B}) on $B^1,\dots,B^N$ ensures that
 \begin{displaymath}
 \mathcal I_N =
 \{i\in\{1,\dots,N\} : B^i\textrm{ is not independent of }\mathcal F_{i - 1}\},
 \end{displaymath}
 where $\mathcal F_0 :=\mathcal F$ and $\mathcal F_i :=\sigma(B^1,\dots,B^i)$ for every $i\in\{1,\dots,N\}$. This is crucial in the proof of Theorem \ref{risk_bound}.
 \item The condition on $|\mathcal I_N|$ in Assumption \ref{sparsity_condition} can be understood as a sparsity type condition on the correlation matrix $R$. Clearly under Assumption \ref{sparsity_condition}, for $N$ large enough,
 \begin{equation}\label{sparsity_condition_1}
 \frac{|\mathcal I_N|^p}{N^{p/2}\log(NT)^{p/2}}
 \leqslant\frac{\mathfrak c_{\ref{sparsity_condition}}(p)}{N^6}
 \quad {\rm with}\quad
 \mathfrak c_{\ref{sparsity_condition}}(p) =
 \mathfrak c_{\ref{sparsity_condition}}^p.
 \end{equation}
 Note that, $p$ is also involved in Assumption \ref{condition_m} through $\mathfrak c_T(p)$. Assumption \ref{sparsity_condition} suggests to take $p$ as large as possible. So, the larger $p$, the larger $\mathcal I_N$, but the smaller the {\it authorized} choices of $m$. In other words, $p$ needs to be chosen large, but not that much.
 \\
 \\
 Note also that for Theorem \ref{risk_bound}, the constraint $p\geqslant 12$ may be lightened into $p\geqslant 4$ and Assumption \ref{sparsity_condition} into
 \begin{displaymath}
 |\mathcal I_N|\leqslant
 \mathfrak c_{\ref{sparsity_condition}}
 N^{\frac{1}{2} -\frac{2}{p}}.
 \end{displaymath}
 However, Theorem \ref{bound_adaptive_estimator} is more demanding.
\end{enumerate}
{\bf Examples:}
\begin{enumerate}
 \item Assume that $N\in q\mathbb N^*$ with $q\in\mathbb N^*$, and that $R$ is the block matrix defined by
 \begin{displaymath}
 R :=
 \begin{pmatrix}
  R_1 & & (0)\\
   & \ddots & \\
  (0) & & R_{\frac{N}{q}}
 \end{pmatrix},
 \end{displaymath}
 where $R_1,\dots,R_{N/q}$ are $N/q$ correlation matrices of size $q\times q$. For instance, if the number of $R_i$'s not equal to $\mathbf I$ is of order lower than $N^{1/2 - 6/p}$, then the matrix $R$ fulfills Assumption \ref{sparsity_condition}. For $q = 2$,
 \begin{displaymath}
 R_i =
 \begin{pmatrix}
  1 & \rho_i\\
  \rho_i & 1
 \end{pmatrix}
 \quad {\rm with}\quad
 \rho_i\in [-1,1]
 \end{displaymath}
 for every $i\in\{1,\dots,N/2\}$. In this special case, $R$ fulfills Assumption \ref{sparsity_condition} if and only if the number of non-zero $\rho_i$'s is of order lower than $N^{1/2 - 6/p}$.
 \item Assume that $N\in 2\mathbb N^*$ and that $R =\mathbf I +\mathbf Q$, where
 \begin{displaymath}
 \mathbf Q :=
 \begin{pmatrix}
  (0) & (0) & Q^*\\
  (0) & (0) & (0)\\
  Q & (0) & (0)
 \end{pmatrix},
 \end{displaymath}
 $Q$ is a correlation matrix of size $r\times r$, and $r\in\{1,\dots,N/2\}$. If $r = r(N)$ is of order lower than $N^{1/2 - 6/p}$, then the matrix $R$ fulfills Assumption \ref{sparsity_condition}.
 \\
 \\
 Note that $R$ is a Toeplitz matrix when
 \begin{displaymath}
 Q =
 \begin{pmatrix}
 0 & 0 & 0 & \cdots & \cdots & 0\\
 q_1 & 0 & 0 & \ddots & & \vdots\\
 q_2 & q_1 & \ddots & \ddots & \ddots & \vdots\\
 \vdots & \ddots & \ddots & \ddots & 0 & 0\\
 \vdots & & \ddots & q_1 & 0 & 0\\
 q_r & \cdots & \cdots & q_2 & q_1 & 0
 \end{pmatrix}
 \quad {\rm with}\quad
 q_1,\dots,q_r\in [-1,1].
 \end{displaymath}
\end{enumerate}
%

% Theorem : Risk bound.

%
\begin{theorem}\label{risk_bound}
Under Assumptions \ref{assumption_basis}, \ref{condition_m} and \ref{sparsity_condition}, there exists a deterministic constant $\mathfrak c_{\ref{risk_bound}} > 0$, not depending on $m$ and $N$, such that
\begin{displaymath}
\mathbb E(\|\widetilde b_m - b_I\|_{N}^{2})
\leqslant
\min_{\tau\in\mathcal S_m}\|\tau - b_I\|_{f_T}^{2} +
\frac{\mathfrak c_{\ref{risk_bound}}m}{NT}
\left(1 +\frac{1}{N}\sum_{i\not= k}|R_{i,k}|\right) +
\frac{\mathfrak c_{\ref{risk_bound}}}{N}.
\end{displaymath}
\end{theorem}
\noindent
As usual for a nonparametric estimator, the risk bound involves a bias term
\begin{displaymath}
\min_{\tau\in\mathcal S_m}\|\tau - b_I\|_{f_T}^{2},
\end{displaymath}
and a variance term of order $m/(NT)$ if
\begin{displaymath}
\frac{1}{N}
\sum_{i,k = 1}^{N}|R_{i,k}| 
\end{displaymath}
is bounded by a constant which does not depend on $N$. The last term is of order $1/N$ and gathers all negligible quantities. The larger $m$, the better the approximation of $b_I$ in $\mathcal S_m$ and the smaller the bias. On the opposite, the variance increases with $m$. This is why a compromise must be done, either theoretically as in Section 2.4 of Comte and Genon-Catalot \cite{CGC20} from which consistency follows, or by a model selection procedure, as described hereafter.

% Section : Model selection.

%
\section{Model selection}\label{model_selection_section}
Throughout this section, $(\varphi_1,\dots,\varphi_{N_T})$ and the $R_{i,k}$'s fulfill the following additional assumptions.
%

% Assumption : Additional assumption on the basis.

%
\begin{assumption}\label{additional_assumption_basis}
The $\varphi_j's$ satisfy the two following (additional) conditions:
\begin{enumerate}
 \item There exists a deterministic constant $\mathfrak c_{\varphi}\geqslant 1$, not depending on $N$, such that for every $m\in\{1,\dots,N_T\}$,
 \begin{displaymath}
 L(m) =
 1\vee\left(\sup_{x\in I}\sum_{j = 1}^{m}\varphi_j(x)^2\right)
 \leqslant\mathfrak c_{\varphi}^{2}m.
 \end{displaymath}
 \item For every $m,m'\in\{1,\dots,N_T\}$, if $m > m'$, then $\mathcal S_{m'}\subset\mathcal S_m$.
\end{enumerate}
\end{assumption}
\noindent
{\bf Remark.} Note that Assumption \ref{additional_assumption_basis}.(2) is fulfilled when
\begin{equation}\label{nested_bases}
\mathcal S_{m + 1} =\mathcal S_m + {\rm span}\{\varphi_{m + 1}\}
\textrm{ $;$ }\forall m\in\{1,\dots,N_T\}.
\end{equation}
For instance, the spaces generated by the trigonometric basis or by the Hermite basis, both defined in Section \ref{section_numerical_experiments}, satisfy (\ref{nested_bases}).
%

% Assumption : Additional assumption on the Brownian montions.

%
\begin{assumption}\label{additional_assumption_Brownians}
There exists a deterministic constant $\mathfrak m_{\ref{additional_assumption_Brownians}} > 0$, not depending on $N$, such that
\begin{displaymath}
\|R\|_{\rm op}\leqslant
\mathfrak m_{\ref{additional_assumption_Brownians}}.
\end{displaymath}
\end{assumption}
\noindent
{\bf Examples (continued).} Since $R$ is a symmetric matrix, there exist an orthogonal matrix $P$ and a diagonal matrix $D$ such that $R = PDP^*$. Then,
\begin{displaymath}
\|R\|_{\rm op} =\|D\|_{\rm op} =\sup_{\lambda\in {\rm sp}(R)}|\lambda|.
\end{displaymath}
So, the matrix $R$ fulfills Assumption \ref{additional_assumption_Brownians} if and only if there exists a constant $\mathfrak m > 0$, not depending on $N$, such that $|\lambda|\leqslant\mathfrak m$ for every $\lambda\in {\rm sp}(R)$. Moreover, note that since $R = PDP^*$,
\begin{displaymath}
\mathcal I_N =\left\{i\in\{2,\dots,N\} :
\exists k\in\{1,\dots,i - 1\}\textrm{ such that }
\sum_{r = 1}^{N}D_{r,r}P_{i,r}P_{k,r}\neq 0\right\}.
\end{displaymath}
\begin{enumerate}
 \item (continued) Assume that $q = 2$. For every $\lambda\in\mathbb R$,
 \begin{displaymath}
 \det(R -\lambda\mathbf I) =
 \prod_{i = 1}^{\frac{N}{2}}\det(R_i -\lambda\mathbf I) =
 \prod_{i = 1}^{\frac{N}{2}}(1 -\lambda -\rho_i)(1 -\lambda +\rho_i),
 \end{displaymath}
 and then
 \begin{displaymath}
 {\rm sp}(R) =
 \left\{1\pm\rho_i\textrm{ $;$ }
 i = 1,\dots,\frac{N}{2}\right\}.
 \end{displaymath}
 So, the matrix $R$ fulfills Assumption \ref{additional_assumption_Brownians}:
 \begin{displaymath}
 \|R\|_{\rm op} =
 \max_{i\in\{1,\dots,N/2\}}|1\pm\rho_i|\leqslant 2.
 \end{displaymath}
 More generally, assume that $q\geqslant 2$. Since
 \begin{displaymath}
 \det(R -\lambda\mathbf I) =
 \prod_{i = 1}^{\frac{N}{q}}\det(R_i -\lambda\mathbf I)
 \textrm{ $;$ }
 \forall\lambda\in\mathbb R,
 \end{displaymath}
 then
 \begin{displaymath}
 {\rm sp}(R) =\bigcup_{i = 1}^{\frac{N}{q}}{\rm sp}(R_i).
 \end{displaymath}
 So, the matrix $R$ fulfills Assumption \ref{additional_assumption_Brownians}:
 \begin{eqnarray*}
  \|R\|_{\rm op} & = &
  \sup_{\lambda\in {\rm sp}(R)}|\lambda| =
  \max_{i\in\{1,\dots,N/q\}}\left\{\sup_{\lambda\in {\rm sp}(R_i)}|\lambda|\right\} =
  \max_{i\in\{1,\dots,N/q\}}\|R_i\|_{\rm op}\\
  & \leqslant &
  \max_{i\in\{1,\dots,N/q\}}\left(\sum_{k,\ell = 1}^{q}[R_i]_{k,\ell}^{2}\right)^{\frac{1}{2}}
  \leqslant q.
 \end{eqnarray*}
 \item (continued) Let us show that $R$ fulfills Assumption \ref{additional_assumption_Brownians} if and only if
 \begin{equation}\label{additional_assumption_Brownians_Toeplitz}
 \sup\{|\lambda|\textrm{ $;$ }\lambda\in\mathbb R\textrm{ such that }
 \det((1 -\lambda)^2\mathbf I_r - Q^*Q) = 0\}
 \textrm{ is bounded by a constant which doesn't depend on $N$.}
 \end{equation}
 For any $\lambda\in\mathbb R$,
 \begin{eqnarray*}
  \det(R -\lambda\mathbf I_N) & = &
  \left|
  \begin{array}{ccc}
   (1 -\lambda)\mathbf I_{\frac{N}{2}} & \overline Q^*\\
   \overline Q & (1 -\lambda)\mathbf I_{\frac{N}{2}}
  \end{array}\right|
  \quad {\rm with}\quad
  \overline Q =
  \begin{pmatrix}
   (0) & (0)\\
   Q & (0)
  \end{pmatrix}
  \in\mathcal M_{\frac{N}{2}}(\mathbb R)\\
  & = &
  \det((1 -\lambda)^2\mathbf I_{\frac{N}{2}} -\overline Q^*\times\overline Q).
 \end{eqnarray*}
 Moreover,
 \begin{displaymath}
 \overline Q^*\times\overline Q =
 \begin{pmatrix}
  Q^*Q & (0)\\
  (0) & (0)
 \end{pmatrix},
 \end{displaymath}
 leading to
 \begin{eqnarray*}
  \det((1 -\lambda)^2\mathbf I_{\frac{N}{2}} -\overline Q^*\times\overline Q)
  & = &
  \left|
  \begin{array}{ccc}
   (1 -\lambda)^2\mathbf I_r - Q^*Q & (0)\\
   (0) & (1 -\lambda)^2\mathbf I_{\frac{N}{2} - r}
  \end{array}\right|\\
  & = &
  (1 -\lambda)^{N - 2r}\det((1 -\lambda)^2\mathbf I_r - Q^*Q).
 \end{eqnarray*}
 So,
 \begin{displaymath}
 \det(R -\lambda\mathbf I_N) =
 (1 -\lambda)^{N - 2r}\det((1 -\lambda)^2\mathbf I_r - Q^*Q),
 \end{displaymath}
 and then $R$ fulfills Assumption \ref{additional_assumption_Brownians} if and only if $Q$ fulfills (\ref{additional_assumption_Brownians_Toeplitz}). Finally, for instance, assume that
 \begin{displaymath}
 Q =
 \begin{pmatrix}
  q_1 & & (0)\\
   & \ddots & \\
  (0) & & q_r
 \end{pmatrix}
 \quad {\rm with}\quad
 q_1,\dots,q_r\in [-1,1].
 \end{displaymath}
 So,
 \begin{eqnarray*}
  \det((1 -\lambda)^2\mathbf I_r - Q^*Q) & = &
  \det((1 -\lambda)^2\mathbf I_r - Q^2)\\
  & = &
  \prod_{i = 1}^{r}(1 -\lambda - q_i)(1 -\lambda + q_i),
 \end{eqnarray*}
 and then $Q$ fulfills (\ref{additional_assumption_Brownians_Toeplitz}) because $q_1,\dots,q_r\in [-1,1]$.
\end{enumerate}
Let us consider
\begin{displaymath}
\widehat m =
\arg\min_{m\in\widehat{\mathcal M}_N}
\{-\|\widehat b_m\|_{N}^{2} + {\rm pen}(m)\},
\end{displaymath}
where
\begin{displaymath}
{\rm pen}(m) :=
\mathfrak c_{\rm cal}
\frac{m}{NT}\left(1 +\frac{1}{N}\sum_{i\neq k}|R_{i,k}|\right)
\textrm{ $;$ }
\forall m\in\{1,\dots,N_T\},
\end{displaymath}
$\mathfrak c_{\rm cal} > 0$ is a deterministic constant to calibrate,
\begin{displaymath}
\widehat{\mathcal M}_N :=
\left\{m\in\{1,\dots,N_T\} :
[\mathfrak c_{\varphi}^{2}m(\|\widehat{\bf\Psi}_{m}^{-1}\|_{\rm op}\vee 1)]^2
\leqslant\mathfrak d_T(p)\frac{NT}{\log(NT)}
\right\}
\end{displaymath}
and
\begin{displaymath}
\mathfrak d_T(p) :=
\frac{1}{512\mathfrak c_{\varphi}^{4}T(1 + p/2)}.
\end{displaymath}
Consider also the theoretical counterpart
\begin{displaymath}
\mathcal M_N :=
\left\{m\in\{1,\dots,N_T\} :
[\mathfrak c_{\varphi}^{2}m(\|\mathbf\Psi_{m}^{-1}\|_{\rm op}\vee 1)]^2
\leqslant\frac{\mathfrak d_T(p)}{4}\cdot\frac{NT}{\log(NT)}
\right\}
\textrm{ of }
\widehat{\mathcal M}_N.
\end{displaymath}
%

% Theorem : Bound on the adaptive least squares projection estimator.

%
\begin{theorem}\label{bound_adaptive_estimator}
Under Assumptions \ref{assumption_basis}, \ref{sparsity_condition}, \ref{additional_assumption_basis} and \ref{additional_assumption_Brownians}, there exist deterministic constants $\kappa_0,\mathfrak c_{\ref{bound_adaptive_estimator}} > 0$, not depending on $N$, such that $\mathfrak c_{\rm cal}\geqslant\kappa_0$ and
\begin{displaymath}
\mathbb E(\|\widehat b_{\widehat m} - b_I\|_{N}^{2})\leqslant
\mathfrak c_{\ref{bound_adaptive_estimator}}\min_{m\in\mathcal M_N}
\left\{\min_{\tau\in\mathcal S_m}
\|\tau - b_I\|_{f_T}^{2} +\frac{m}{NT}\left(1 +\frac{1}{N}\sum_{i\neq k}|R_{i,k}|\right)
\right\} +
\frac{\mathfrak c_{\ref{bound_adaptive_estimator}}}{N}.
\end{displaymath}
\end{theorem}
\noindent
It follows from Theorem \ref{bound_adaptive_estimator} that the adaptive estimator $\widehat b_{\widehat m}$ automatically reaches a squared bias-variance compromise on the collection $\mathcal M_N$.
\\
\\
{\bf Remark.} Note that the constant $\kappa_0$ is given at the end of the proof of Lemma \ref{bound_empirical_process}.
%

% Section : Numerical experiments.

%
\section{Numerical experiments}\label{section_numerical_experiments}
In this section, we study the influence of dependency on the performance of the adaptive estimator. We consider two bases:
\begin{itemize}
 \item The cosine basis on $I = [a,b]$, defined by $\varphi_1(x) := (b - a)^{-1/2}\mathbf 1_{[a,b]}(x)$, $\varphi_j(x) := (2/(b - a))^{1/2}\cos(\pi j(x - a)/(b - a))\mathbf 1_{[a,b]}(x)$ for $j\geqslant 2$. The interval $[a,b]$ is chosen different for each model. The basis is orthonormal and fulfills $\sum_{j = 1}^{m}\varphi_{j}^{2}(x)\leqslant 2m$.
 \item The Hermite basis on $I =\mathbb R$, defined from the Hermite polynomials $H_j$ and given by
 \begin{displaymath}
 H_j(x) := (-1)^je^{x^2}\frac{d^j}{dx^j}(e^{-x^2}),\quad
 \varphi_j(x) := c_{j - 1}H_{j - 1}(x)e^{-x^2/2},\quad
 c_j =\left(2^j j!\sqrt\pi\right)^{-1/2}.
 \end{displaymath}
 The sequence $(\varphi_j)_{j\geqslant 0}$ is an orthonormal bounded basis of $\mathbb L^2(\mathbb R,dx)$ with $|\varphi_j(x)|\leqslant 1/\pi^{1/4}$ (see Indritz \cite{INDRITZ61}). It is proved in Comte and Lacour \cite{CL21} (see Lemma 1) that $L(m)\leqslant K\sqrt m$ for some constant $K$.
\end{itemize}
We experiment five models, where $I$ is the chosen domain of representation for the Hermite basis, and the basis support for the cosine basis:
\begin{enumerate}
 \item Hyperbolic diffusion, $b_1(x) = -\theta x$ and $\sigma_1(x) =\gamma\sqrt{1 + x^2}$, with $\theta = 2$ and $\gamma =\sqrt{1/2}$, $I_1 = [-0.9,0.8]$.
 \item Hyperbolic tangent of an Ornstein-Uhlenbeck process,
 \begin{displaymath}
 b_2(x) = (1 - x^2)\left(-\frac{r}{2}{\rm atanh}(x) -
 \frac{\gamma^2}{4}x\right),\quad
 \sigma_2(x) =\frac{\gamma}{2}(1 - x^2),\quad r = 2,\quad\gamma = 2,
 \quad I_2 = [-0.9,0.9].
 \end{displaymath}
 \item Exponential of an Ornstein-Ulhenbeck process,
 \begin{displaymath}
 b_3(x) = x\left(-\frac{r}{2}\log(x^+) +\frac{\gamma^2}{8}\right),\quad
 \sigma_3(x) =\frac{\gamma}{2}x^+,\quad
 r = 1,\quad\gamma = 2,\quad I_3 = [0.44,2].
 \end{displaymath}
 \item $X_t = G_1(\xi_t)$ with $d\xi_t =\alpha(\xi_t)dt + dW_t$, $\alpha(x) = -\theta x/\sqrt{1 + c^2x^2}$, $G_1(x) = {\rm asinh}(c\xi_t)$, $\theta = 3$, $c = 2$, $I_4 = [-1.15,1.15]$, leading to
 \begin{displaymath}
 b_4(x) = -\left(\theta +\frac{c^2}{2\cosh(x)}\right)\frac{\sinh(x)}{\cosh^2(x)},\quad
 \sigma_4(x) =\frac{c}{\cosh(x)}.
 \end{displaymath}
 \item $X_t = G_2(\xi_t)$ with $\xi$ as previously and $G_2(x) = {\rm asinh}(x - 5) + {\rm asinh}(x + 5)$, leading to
 \begin{displaymath}
 b_5(x) = G_2'(H(x))\alpha(H(x)) +\frac{1}{2}G_2''(H(x)),\quad
 \sigma_5(x) =\frac{1}{\sqrt{1 + (H(x) - 5)^2}} +\frac{1}{\sqrt{1 + (H(x) + 5)^2}},
 \end{displaymath}
 where $\theta = 1$, $c = 10$ (in the definition of $\alpha$), $I_5 = [-4,4]$ and
 \begin{displaymath}
 H(x) = G_{2}^{-1}(x) =\frac{1}{\sqrt{2}\sinh(x)}
 \left[(49 +\cosh(x))\sinh(x)^2 + 100(1 -\cosh(x)\right]^{1/2}.
 \end{displaymath}
\end{enumerate}
Models 1, 4, 5 are simulated by Euler scheme with step $\Delta$, directly for $X$ in example 1 or for $\xi$ in examples 4 and 5, with transformations $G_1$ and $G_2$ in a second stage. The underlying Ornstein-Uhlenbeck processes in models 2 and 3 are generated by exact autoregressive scheme with step $\Delta$. Details can be found in Comte and Genon-Catalot \cite{CGC21} for examples 1, 2, 3 and in Comte al al. \cite{CGCR07} for examples 4 and 5.
\begin{figure}
\includegraphics[width=16cm, height=6cm]{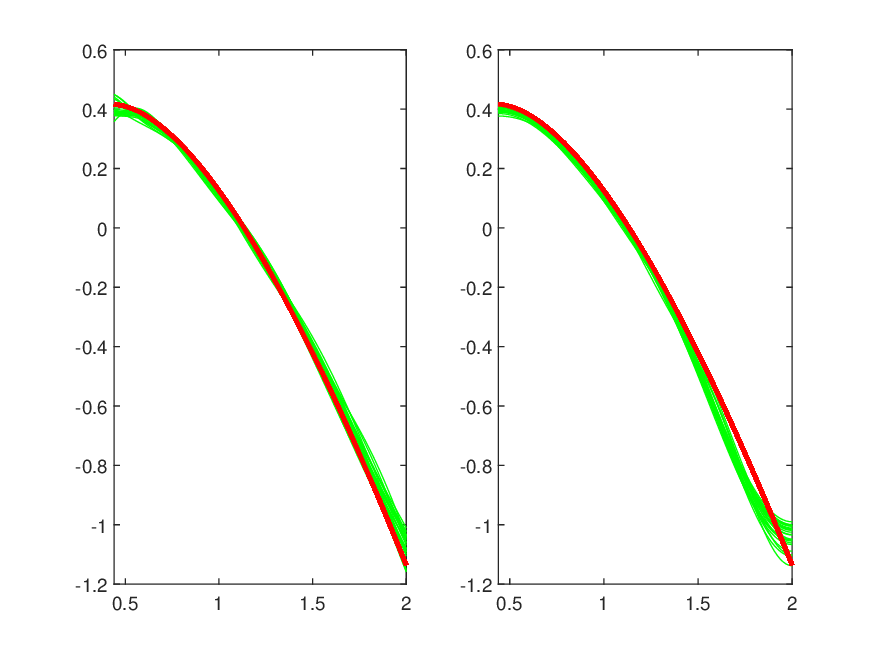}
\caption{Example 3. True functions in bold red and beam of 25 estimated drift $b_3$ with Hermite (left) and cosine (right) bases, $\rho = 0.5$. The MISE$\times 100$ are $0.12$, $0.33$ and the mean of selected dimensions are $8.4$, $4.3$.}
\label{fig1}
\end{figure}
\begin{figure}
\includegraphics[width=16cm, height=6cm]{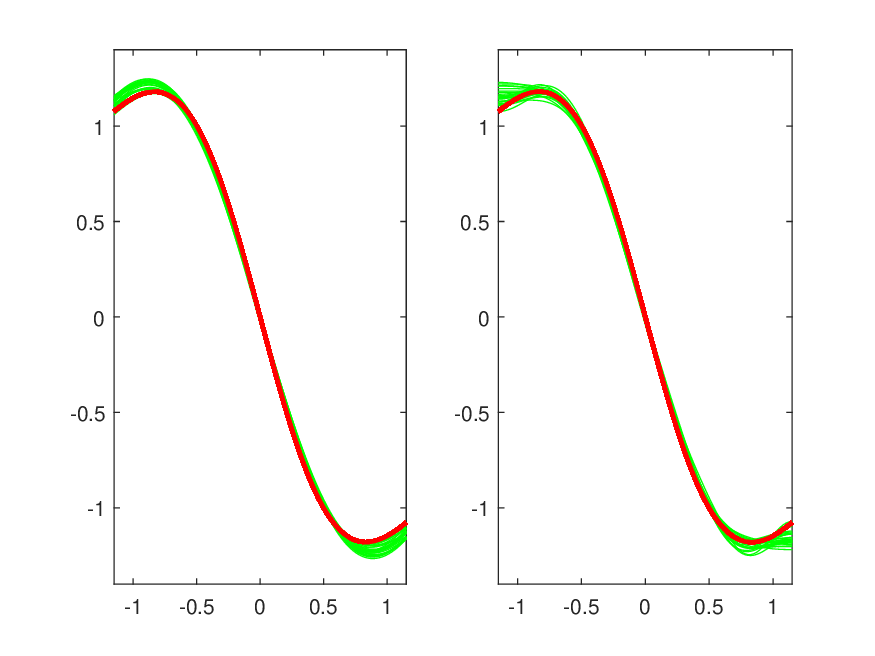}
\caption{Example 4. True functions in bold red and beam of 25 estimated drift $b_4$ with Hermite (left) and cosine (right) bases, $\rho = 0.5$. The MISE$\times 100$ are $0.36$, $0.26$ and the mean of selected dimensions are $4.4$, $6.0$.}
\label{fig2}
\end{figure}
\begin{figure}
\begin{tabular}{cc}
\includegraphics[scale=0.5]{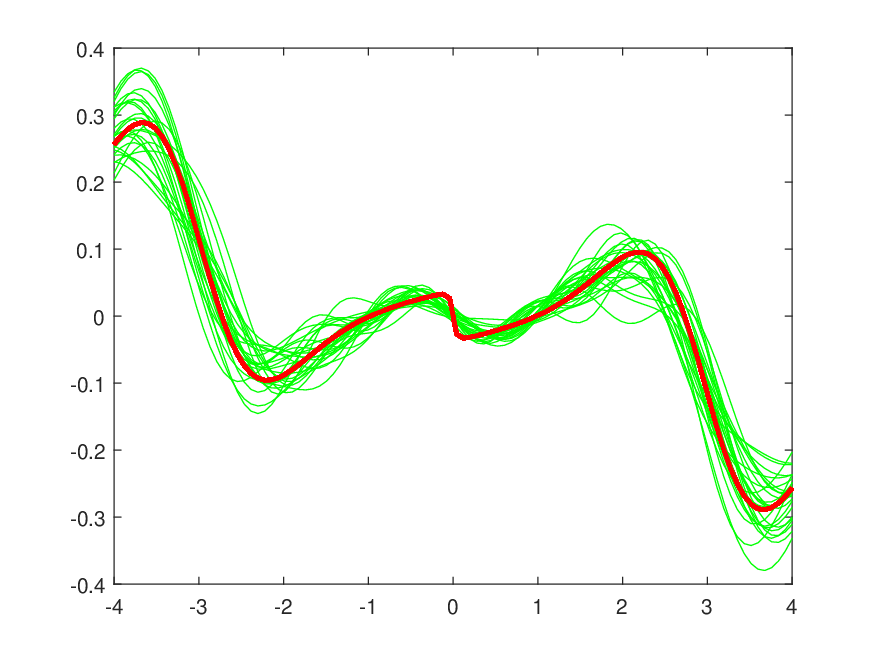} &
\includegraphics[scale=0.5]{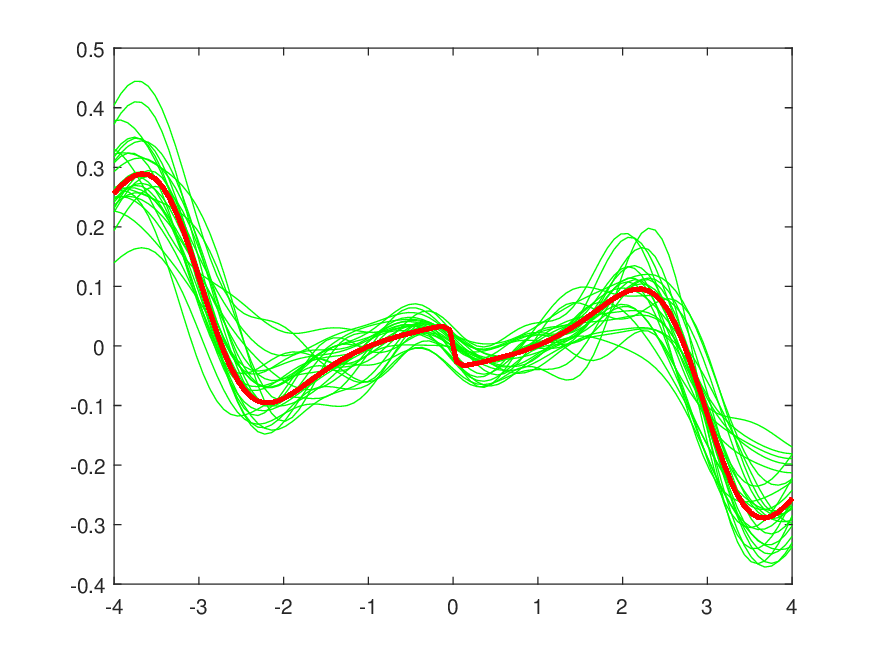} 
\end{tabular}
\caption{Example 5. True functions in bold red and beam of 25 estimated drift $b_5$ with Hermite basis for $\rho = 0$ (left) and $\rho = 0.9$ (right). The MISE$\times 100$ are $0.87$, $1.67$ and the mean of selected dimensions are $12.5$ in both cases.}
\label{fig3}
\end{figure}
The dependency is contained in the Toeplitz variance matrix\footnote{This matrix is indeed a correlation matrix as it is the variance of a stationary AR(1) process, $X_t =\rho X_{t - 1} +\varepsilon_t$, for i.i.d. centered $\varepsilon_t$'s.} $R(\rho) := (\rho^{|i - j|})_{1\leqslant i,j\leqslant N}$ for different values of $\rho$. The choice $\rho = 0$ corresponds to the independent case, and we also experiment $\rho = 0.5$ (mild dependency) and $\rho = 0.9$ (strong correlations). Assumption \ref{sparsity_condition} is not fulfilled but we can consider that the coefficients are in fact null when $|i - j|$ is large enough. The orders of some quantities related to $R(\rho)$ are given in Table \ref{tab0}, and clearly, $\sum_{i,k}R(\rho)_{i,k}/N$ and $\|R(\rho)\|_{\rm op}$ are very close.
\begin{table}[h!]
\begin{center}
\begin{tabular}{c|ccc} 
 & $\rho = 0$ & $\rho = 0.5$ & $\rho = 0.9$\\
 \hline
 $N^{-1}\sum_{i,j}R_{i,j}(\rho)$ & 1 & 2.96 & 17.2\\
 $\|R(\rho)\|_{\rm op}$ & 1 & 2.99 & 17.9\\
\multicolumn{4}{c}{$\;$}
\end{tabular}
\caption{Order of different quantities associated to the matrix $R(\rho)$.}\label{tab0}
\end{center}
\end{table}
The penalty term is computed as in Comte and Genon-Catalot \cite{CGC20}, by an empirical version which directly takes dependence into account without requiring any information on $R$:
\begin{displaymath}
\widehat{\rm pen}(m) :=
\kappa\frac{m}{NT}\|\widehat{\bf\Psi}_{m}^{-1}
\widehat{\bf\Psi}_{m,\sigma}\|_{\rm op}
\quad {\rm with}\quad
\widehat{\bf\Psi}_{m,\sigma} =
(\langle\sigma\varphi_j,\sigma\varphi_{\ell}\rangle_N)_{j,\ell}
\end{displaymath}
and $\kappa = 2$ for both bases. Then, the model $\widehat m$ is chosen as the minimizer of $-\|\widehat b_m\|_{N}^{2} +\widehat{\rm pen}(m)$ for $m\leqslant 10$ (resp. $m\leqslant 20$) for the Hermite basis, except in example 5, where we set $m\leqslant 15$ because otherwise the selected dimension was systematically the maximal one (resp. for the cosine basis), such that
\begin{displaymath}
m\|\widehat{\bf\Psi}_{m}^{-1}\|_{\rm op}^{1/4}\leqslant NT
\textrm{ (empirical collection of models)}.
\end{displaymath}
We present illustrations of the estimation procedure obtained from simulated paths in Figures \ref{fig1}, \ref{fig2}, \ref{fig3} for examples 3, 4 and 5. Figures \ref{fig1} and \ref{fig2} allow the comparison of the results obtained for Hermite (left pictures) and cosine (right pictures) bases, for $\rho = 0.5$. Figure \ref{fig3}  shows the difference of estimation in the Hermite basis when $\rho = 0$ (left picture) and when $\rho = 0.9$ (right picture). The scenario is the same in the three figures: $N = 200$ and $T=100$ (with $1000$ observations with step $\Delta = 0.1$ for each path). The MISE over the 25 repetitions are given, together with the mean of the selected dimensions. We can see that the examples are quite different, and that the estimation method works in a convincing way, even for strong dependency ($\rho = 0.9$).\\
We also illustrate on the scenario $N = 100$ and $T = 100$ (with $1000$ observations with step $\Delta = 0.1$) for each path, which was a middle scenario in Comte and Genon-Catalot \cite{CGC20}, the influence of the value of $\rho$ on the MISE computed over 200 repetitions: the results are given in Table \ref{tab1}. We see that the MISE increases when $\rho$ increases, slightly from $\rho = 0$ to $\rho = 0.5$ and much more importantly from $\rho = 0.5$ to $\rho = 0.9$. On the contrary, the selected dimensions for each basis are rather unchanged in these different cases. This suggests that bias and variance increase simultaneously and proportionally. The Hermite basis gives lower MISEs for examples 1 to 3, and the cosine basis wins for examples 4 and 5.
\begin{table}[h!]
\begin{center}
\begin{tabular}{cc|cc|cc|cc}
 && \multicolumn{2}{c}{$\rho = 0$} & \multicolumn{2}{c}{$\rho = 0.5$} & \multicolumn{2}{c}{$\rho = 0.9$}\\
 Ex. && Hermite & Cosine & Hermite & Cosine & Hermite & Cosine\\
 \hline
 Ex.1& MISE & $0.11_{(0.08)}$ & $0.80_{(0.24)}$ & $0.13_{(0.10)}$ & $0.83_{(0.29)}$ & $0.63_{(0.54)}$ & $1.38_{(6.29)}$\\
 & Dim  & $6.2_{(0.8)}$ & $6.2_{(1.5)}$ & $6.1_{(0.5)}$ & $6.2_{(1.5)}$ & $6.3_{(1.1)}$ & $6.3_{(1.2)}$\\
 \hline 
 Ex.2& MISE & $0.78_{(0.18)}$ & $0.95_{(0.19)}$ & $0.78_{(0.18)}$ & $0.94_{(0.18)}$ & $1.02_{(0.47)}$ & $1.18_{(0.49)}$\\
 & Dim  & $6.1_{(0.5)}$ & $10.3_{(1.9)}$ & $6.1_{(0.5)}$ & $10.4_{(2.0)}$ & $6.1_{(0.5)}$ & $10.3_{(2.0)}$\\
 \hline 
 Ex.3& MISE  & $0.22_{(0.16)}$ & $0.34_{(0.11)}$ & $0.21_{(0.16)}$ & $0.37_{(0.14)}$ & $0.44_{(0.41)}$ & $0.55_{(0.43)}$\\
 & Dim  & $7.8_{(0.7)}$ & $4.1_{(0.4)}$ & $7.7_{(0.7)}$ & $4.1_{(0.4)}$ & $7.8_{(0.8)}$ & $4.2_{(0.6)}$\\
 \hline 
 Ex.4& MISE & $0.41_{(0.18)}$ & $0.35_{(0.15)}$ & $0.46_{(0.21)}$ & $0.39_{(0.18)}$ & $1.06_{(0.66)}$ & $1.02_{(0.60)}$\\
 & Dim  & $4.4_{(0.8)}$ & $5.1_{(1.4)}$ & $4.4_{(0.8)}$ & $5.2_{(1.3)}$ & $4.9_{(1.3)}$ & $5.6_{(1.5)}$\\
 \hline 
 Ex.5& MISE & $1.55_{(0.85)}$ & $1.49_{(0.61)}$ & $1.81_{(0.97)}$ & $1.58_{(0.62)}$ & $4.14_{(4.22)}$ & $3.15_{(3.55)}$\\
 & Dim  & $11.2_{(1.3)}$ & $6.3_{(0.8)}$ & $11.4_{(1.4)}$ & $6.2_{(0.6)}$ & $11.5_{(1.5)}$ & $6.3_{(0.8)}$\\
 \multicolumn{8}{c}{$\;$}
\end{tabular}
\caption{100 MISE (with 100 Std in parenthesis) and mean selected dimensions (with StD in parenthesis) for the examples 1 to 5, $N = 100$ and $T = 100$, for Hermite and cosine bases and 3 values of $\rho$ (0 for independence, $\rho = 0.9$ for strong dependency).}\label{tab1}
\end{center}
\end{table}
%

% Section : Acknowledgments.

%
\section*{Acknowledgments}
This research has been conducted as part of the MALNOS project funded by Labex MME-DII (ANR11-LBX-0023- 01).
%

% Section : Proofs

%
\section{Proofs}\label{proofs_section}
%

% Subsection : Proof of Lemma preliminary_lemma_variance.

%
\subsection{Proof of Lemma \ref{preliminary_lemma_variance}}
First of all, let us show that the symmetric matrix $\mathbf\Psi_{m,\sigma}$ is nonnegative. Indeed, for every $y\in\mathbb R^m$,
\begin{eqnarray*}
 y^*\mathbf\Psi_{m,\sigma}y & = &
 \frac{1}{NT}\sum_{j,\ell = 1}^{m}y_jy_{\ell}
 \sum_{i,k = 1}^{N}\mathbb E\left(
 \left(\int_{0}^{T}\sigma(X_{s}^{i})\varphi_j(X_{s}^{i})dB_{s}^{i}\right)
 \left(\int_{0}^{T}\sigma(X_{s}^{k})\varphi_{\ell}(X_{s}^{k})dB_{s}^{k}\right)\right)\\
 & = &
 \frac{1}{NT}\mathbb E\left[\left(
 \sum_{i = 1}^{N}\int_{0}^{T}\sigma(X_{s}^{i})\tau_y(X_{s}^{i})dB_{s}^{i}\right)^2\right]
 \geqslant 0
 \quad {\rm with}\quad
 \tau_y(.) :=\sum_{j = 1}^{m}y_j\varphi_j(.).
\end{eqnarray*}
On the one hand, since $\mathbf\Psi_{m,\sigma}$ is a nonnegative matrix, since $d\langle B^i,B^k\rangle_t = R_{i,k}dt$ for every $i,k\in\{1,\dots,N\}$, and by the stochastic integration by parts formula,
\begin{eqnarray*}
 {\rm trace}(\mathbf\Psi_{m}^{-1}\mathbf\Psi_{m,\sigma})
 & \leqslant &
 \|\mathbf\Psi_{m}^{-1}\|_{\rm op}
 {\rm trace}(\mathbf\Psi_{m,\sigma}) =
 \frac{1}{NT}
 \|\mathbf\Psi_{m}^{-1}\|_{\rm op}
 \sum_{j = 1}^{m}\mathbb E\left[\left(\sum_{i = 1}^{N}\int_{0}^{T}
 \sigma(X_{s}^{i})\varphi_j(X_{s}^{i})dB_{s}^{i}\right)^2\right]\\
 & = &
 \frac{1}{NT}\|\mathbf\Psi_{m}^{-1}\|_{\rm op}\sum_{j = 1}^{m}\sum_{i,k = 1}^{N}
 \int_{0}^{T}R_{i,k}
 \mathbb E(\sigma(X_{s}^{i})\varphi_j(X_{s}^{i})
 \sigma(X_{s}^{k})\varphi_j(X_{s}^{k}))ds\\
 & \leqslant &
 \frac{1}{N}\|\mathbf\Psi_{m}^{-1}\|_{\rm op}
 \left(N +\sum_{i\neq k}|R_{i,k}|\right)
 \sum_{j = 1}^{m}\int_{-\infty}^{\infty}
 \sigma(x)^2\varphi_j(x)^2f_T(x)dx\\
 & \leqslant &
 \mathfrak c_1\|\mathbf\Psi_{m}^{-1}\|_{\rm op}L(m)
 \left(1 +\frac{1}{N}\sum_{i\neq k}|R_{i,k}|\right)
 \quad {\rm with}\quad
 \mathfrak c_1 =\int_{-\infty}^{\infty}\sigma(x)^2f_T(x)dx.
\end{eqnarray*}
On the other hand, assume now that $\sigma$ is bounded. Again, since $d\langle B^i,B^k\rangle_t = R_{i,k}dt$ for every $i,k\in\{1,\dots,N\}$, and by the stochastic integration by parts formula, for every $y\in\mathbb R^m$,
\begin{eqnarray}
 y^*\mathbf\Psi_{m,\sigma}y & = &
 \frac{1}{NT}\mathbb E\left[
 \left(\sum_{i = 1}^{N}\int_{0}^{T}
 \sigma(X_{s}^{i})\tau_y(X_{s}^{i})dB_{s}^{i}\right)^2\right]
 \nonumber\\
 & \leqslant &
 \frac{1}{T}\left(1 +\frac{1}{N}\sum_{i\neq k}|R_{i,k}|\right)
 \int_{0}^{T}\mathbb E(\sigma(X_s)^2\tau_y(X_s)^2)ds
 \nonumber\\
 & \leqslant &
 \left(1 +\frac{1}{N}\sum_{i\neq k}|R_{i,k}|\right)
 \int_{-\infty}^{\infty}\sigma(x)^2
 \left(\sum_{j = 1}^{m}y_j\varphi_j(x)\right)^2f_T(x)dx
 \nonumber\\
 \label{preliminary_lemma_variance_1}
 & \leqslant &
 \|\sigma\|_{\infty}^{2}
 \left(1 +\frac{1}{N}\sum_{i\neq k}|R_{i,k}|\right)
 \|\mathbf\Psi_{m}^{1/2}y\|_{2,m}^{2}.
\end{eqnarray}
Thus, since $\mathbf\Psi_{m,\sigma}$ is nonnegative, and by Inequality (\ref{preliminary_lemma_variance_1}),
\begin{eqnarray*}
 {\rm trace}(\mathbf\Psi_{m}^{-1/2}\mathbf\Psi_{m,\sigma}\mathbf\Psi_{m}^{-1/2})
 & \leqslant &
 m\|\mathbf\Psi_{m}^{-1/2}\mathbf\Psi_{m,\sigma}\mathbf\Psi_{m}^{-1/2}\|_{\rm op}\\
 & &
 \hspace{2cm} =
 m\cdot
 \sup\{y^*\mathbf\Psi_{m,\sigma}y\textrm{ $;$ }
 y\in\mathbb R^m\textrm{ and }\|\mathbf\Psi_{m}^{1/2}y\|_{2,m} = 1\}\\
 & \leqslant &
 m\|\sigma\|_{\infty}^{2}
 \left(1 +\frac{1}{N}\sum_{i\neq k}|R_{i,k}|\right).
\end{eqnarray*}
%

% Subsection : Proof of Theorem risk_bound.

%
\subsection{Proof of Theorem \ref{risk_bound}}
The proof of Theorem \ref{risk_bound} relies on the two following lemmas.
%

% Lemma : Bound on the variance of the remainder term.

%
\begin{lemma}\label{bound_variance_remainder}
There exists a deterministic constant $\mathfrak c_{\ref{bound_variance_remainder}} > 0$, not depending on $m$ and $N$, such that
\begin{displaymath}
\mathbb E(|\widehat{\bf E}_{m}^{*}\widehat{\bf E}_m|^2)
\leqslant
\mathfrak c_{\ref{bound_variance_remainder}}
\frac{mL(m)^2}{N^2}\left[1 +
\left(\frac{1}{N}\sum_{i\not= k}|R_{i,k}|\right)^2\right].
\end{displaymath}
\end{lemma}
%

% Lemma : Deviation probabilities.

%
\begin{lemma}\label{deviation_probabilities}
Consider the event
\begin{displaymath}
\Omega_m :=
\left\{\sup_{\tau\in\mathcal S_m}
\left|\frac{\|\tau\|_{N}^{2}}{\|\tau\|_{f_T}^{2}} - 1\right|\leqslant\frac{1}{2}
\right\}.
\end{displaymath}
Under Assumptions \ref{assumption_basis}, \ref{condition_m} and \ref{sparsity_condition}, there exists a deterministic constant $\mathfrak c_{\ref{deviation_probabilities}} > 0$, not depending on $m$ and $N$, such that
\begin{displaymath}
\mathbb P(\Omega_{m}^{c})
\leqslant
\frac{\mathfrak c_{\ref{deviation_probabilities}}}{N^6}
\quad\textrm{and}\quad
\mathbb P(\widehat\Lambda_{m}^{c})
\leqslant
\frac{\mathfrak c_{\ref{deviation_probabilities}}}{N^6}.
\end{displaymath}
\end{lemma}
%

% Subsubsection : Steps of the proof.

%
\subsubsection{Steps of the proof}
First of all,
\begin{eqnarray}
 \|\widetilde b_m - b_I\|_{N}^{2} & = &
 \|b_I\|_{N}^{2}\mathbf 1_{\widehat\Lambda_{m}^{c}} +
 \|\widehat b_m - b_I\|_{N}^{2}\mathbf 1_{\widehat\Lambda_m}
 \nonumber\\
 \label{risk_bound_1}
 & = &
 U_1 + U_2 + U_3
\end{eqnarray}
where $U_1 :=\|b_I\|_{N}^{2}\mathbf 1_{\widehat\Lambda_{m}^{c}}$,
\begin{displaymath}
U_2 :=
\|\widehat b_m - b_I\|_{N}^{2}\mathbf 1_{\widehat\Lambda_m\cap\Omega_m}
\quad {\rm and}\quad
U_3 :=
\|\widehat b_m - b_I\|_{N}^{2}\mathbf 1_{\widehat\Lambda_m\cap\Omega_{m}^{c}}.
\end{displaymath}
Let us find suitable bounds on $\mathbb E(U_1)$, $\mathbb E(U_2)$ and $\mathbb E(U_3)$.
\begin{itemize}
 \item {\bf Bound on $\mathbb E(U_1)$.} By Cauchy-Schwarz's inequality,
 \begin{eqnarray*}
  \mathbb E(U_1) & \leqslant &
  \mathbb E(\|b_I\|_{N}^{4})^{1/2}\mathbb P(\widehat\Lambda_{m}^{c})^{1/2}
  \leqslant
  \mathbb E\left(\frac{1}{T}\int_{0}^{T}b_I(X_t)^4dt\right)^{1/2}
  \mathbb P(\widehat\Lambda_{m}^{c})^{1/2}\\
  & \leqslant &
  \mathfrak c_1\mathbb P(\widehat\Lambda_{m}^{c})^{1/2} <\infty
  \quad {\rm with}\quad
  \mathfrak c_1 =\left(\int_{-\infty}^{\infty}b_I(x)^4f_T(x)dx\right)^{1/2} <\infty.
 \end{eqnarray*}
 \item {\bf Bound on $\mathbb E(U_2)$.} Let $\Pi_{N,m}(.)$ be the orthogonal projection from $\mathbb L^2(I,f_T(x)dx)$ onto $\mathcal S_m$ with respect to the empirical scalar product $\langle .,.\rangle_N$. Then,
 \begin{equation}\label{risk_bound_2}
 \|\widehat b_m - b_I\|_{N}^{2} =
 \|\widehat b_m -\Pi_{N,m}(b_I)\|_{N}^{2} +
 \min_{\tau\in\mathcal S_m}\|b_I -\tau\|_{N}^{2}.
 \end{equation}
 As in the proof of Comte and Genon-Catalot \cite{CGC20}, Proposition 2.1, on $\Omega_m$,
 \begin{displaymath}
 \|\widehat b_m -\Pi_{N,m}(b_I)\|_{N}^{2} =
 \widehat{\bf E}_{m}^{*}\widehat{\bf\Psi}_{m}^{-1}\widehat{\bf E}_m
 \leqslant 2\widehat{\bf E}_{m}^{*}\mathbf\Psi_{m}^{-1}\widehat{\bf E}_m.
 \end{displaymath}
 So,
 \begin{eqnarray*}
  \mathbb E(\|\widehat b_m -\Pi_{N,m}(b_I)\|_{N}^{2}
  \mathbf 1_{\widehat\Lambda_m\cap\Omega_m}) & \leqslant &
  2\mathbb E\left(\sum_{j,\ell = 1}^{m}[\widehat{\bf E}_m]_j[\widehat{\bf E}_m]_{\ell}
  \mathbf\Psi_{m}^{-1}(j,\ell)\right)\\
  & = &
  \frac{2}{NT}\sum_{j,\ell = 1}^{m}\mathbf\Psi_{m,\sigma}(j,\ell)
  \mathbf\Psi_{m}^{-1}(j,\ell) =
  \frac{2}{NT}{\rm trace}(\mathbf\Psi_{m}^{-1}\mathbf\Psi_{m,\sigma})\\
  & = &
    \frac{2}{NT}{\rm trace}(\mathbf\Psi_{m}^{-1/2}
  \mathbf\Psi_{m,\sigma}
  \mathbf\Psi_{m}^{-1/2}).
 \end{eqnarray*}
 Then, by Equality (\ref{risk_bound_2}) and Lemma \ref{preliminary_lemma_variance},
 \begin{eqnarray*}
  \mathbb E(U_2) & \leqslant &
  \mathbb E\left(
  \min_{\tau\in\mathcal S_m}\|b_I -\tau\|_{N}^{2}\right) +
  \frac{2}{NT}{\rm trace}(\mathbf\Psi_{m}^{-1/2}
  \mathbf\Psi_{m,\sigma}
  \mathbf\Psi_{m}^{-1/2})\\
  & \leqslant &
  \min_{\tau\in\mathcal S_m}\|b_I -\tau\|_{f_T}^{2} +
  \frac{2m}{NT}\|\sigma\|_{\infty}^{2}
  \left(1 +\frac{1}{N}\sum_{i\neq k}|R_{i,k}|\right).
 \end{eqnarray*}
 \item {\bf Bound on $\mathbb E(U_3)$.} By the definition of the event 
 $\widehat\Lambda_m$ and by Lemma \ref{bound_variance_remainder},
 \begin{eqnarray*}
  \mathbb E(\|\widehat b_m -\Pi_{N,m}(b_I)\|_{N}^{2}
  \mathbf 1_{\widehat\Lambda_m\cap\Omega_{m}^{c}})
  & = &
  \mathbb E[(\widehat{\bf E}_{m}^{*}\widehat{\bf\Psi}_{m}^{-1}
  \widehat{\bf E}_m)\mathbf 1_{\widehat\Lambda_m\cap\Omega_{m}^{c}}]\\
  & &
  \hspace{2cm}
  \leqslant
  \mathbb E(\|\widehat{\bf\Psi}_{m}^{-1}\|_{\rm op}
  |\widehat{\bf E}_{m}^{*}\widehat{\bf E}_m|
  \mathbf 1_{\widehat\Lambda_m\cap\Omega_{m}^{c}})\\
  & \leqslant &
  \frac{\mathfrak c_T(p)}{L(m)}\cdot
  \frac{NT}{\log(NT)}
  \mathbb E(|\widehat{\bf E}_{m}^{*}\widehat{\bf E}_m|^2)^{1/2}
  \mathbb P(\Omega_{m}^{c})^{1/2}\\
  & \leqslant &
  \frac{\mathfrak c_2m^{1/2}}{\log(NT)}
  \left(1 +\frac{1}{N}\sum_{i\neq k}|R_{i,k}|\right)\mathbb P(\Omega_{m}^{c})^{1/2}
 \end{eqnarray*}
 where $\mathfrak c_2 > 0$ is a deterministic constant not depending on $m$ and $N$. Then,
 \begin{eqnarray*}
  \mathbb E(U_3) & \leqslant &
  \mathbb E(\|\widehat b_m -\Pi_{N,m}(b_I)\|_{N}^{2}
  \mathbf 1_{\widehat\Lambda_m\cap\Omega_{m}^{c}}) +
  \mathbb E(\|b_I\|_{N}^{2}\mathbf 1_{\widehat\Lambda_m\cap\Omega_{m}^{c}})\\
  & \leqslant &
  \frac{\mathfrak c_2m^{1/2}}{\log(NT)}
  \left(1 +\frac{1}{N}\sum_{i\neq k}|R_{i,k}|\right)\mathbb P(\Omega_{m}^{c})^{1/2} +
  \mathfrak c_1\mathbb P(\Omega_{m}^{c})^{1/2}.
 \end{eqnarray*}
\end{itemize}
So,
\begin{eqnarray*}
 \mathbb E(\|\widetilde b_m - b_I\|_{N}^{2})
 & \leqslant &
 \min_{\tau\in\mathcal S_m}\|b_I -\tau\|_{f_T}^{2}\\
 & &
 \hspace{0.5cm} +
 \left[\frac{2m}{NT}\|\sigma\|_{\infty}^{2} +
 \mathfrak c_2
 \frac{\sqrt{m\mathbb P(\Omega_{m}^{c})}}{\log(NT)}\right]
 \left(1 +\frac{1}{N}\sum_{i\not= k}|R_{i,k}|\right) +
 \mathfrak c_1(
 \mathbb P(\widehat\Lambda_{m}^{c})^{1/2} +\mathbb P(\Omega_{m}^{c})^{1/2}).
\end{eqnarray*}
Therefore, by Lemma \ref{deviation_probabilities}, there exists a deterministic constant $\mathfrak c_3 > 0$, not depending on $m$ and $N$, such that
\begin{displaymath}
\mathbb E(\|\widetilde b_m - b_I\|_{N}^{2})
\leqslant
\min_{\tau\in\mathcal S_m}\|b_I -\tau\|_{f_T}^{2} +
\frac{\mathfrak c_3m}{NT}
\left(1 +\frac{1}{N}\sum_{i\not= k}|R_{i,k}|\right) +
\frac{\mathfrak c_3}{N}.
\end{displaymath}
%

% Subsubsection : Proof of Lemma bound_variance_remainder.

%
\subsubsection{Proof of Lemma \ref{bound_variance_remainder}}
By Jensen's inequality, by Burkholder-Davis-Gundy's inequality, and since $d\langle B^i,B^k\rangle_t = R_{i,k}dt$ for every $i,k\in\{1,\dots,N\}$, there exists a deterministic constant $\mathfrak c_1 > 0$, not depending on $m$ and $N$, such that
\begin{eqnarray*}
 \mathbb E(|\widehat{\bf E}_{m}^{*}\widehat{\bf E}_m|^2) & \leqslant &
 m\sum_{j = 1}^{m}\mathbb E(\widehat{\bf E}_m(j)^4)
 \leqslant
 \frac{\mathfrak c_1m}{N^4T^4}\sum_{j = 1}^{m}\mathbb E\left(
 \left\langle\sum_{i = 1}^{N}\int_{0}^{.}
 \sigma(X_{s}^{i})\varphi_j(X_{s}^{i})dB_{s}^{i}\right\rangle_{T}^{2}\right)\\
 & \leqslant &
 \frac{2\mathfrak c_1m}{N^4T^4}\sum_{j = 1}^{m}
 (\mathbb E(D_{j}^{2}) +\mathbb E(A_{j}^{2})),
\end{eqnarray*}
where
\begin{displaymath}
D_j :=\sum_{i = 1}^{N}\int_{0}^{T}
\sigma(X_{s}^{i})^2\varphi_j(X_{s}^{i})^2ds
\quad {\rm and}\quad
A_j :=\sum_{i\neq k}R_{i,k}\int_{0}^{T}
\sigma(X_{s}^{i})\varphi_j(X_{s}^{i})
\sigma(X_{s}^{k})\varphi_j(X_{s}^{k})ds
\end{displaymath}
for every $j\in\{1,\dots,m\}$. On the one hand, by Jensen's inequality,
\begin{eqnarray*}
 \sum_{j = 1}^{m}\mathbb E(D_{j}^{2})
 & \leqslant &
 NT\sum_{j = 1}^{m}\sum_{i = 1}^{N}\int_{0}^{T}
 \mathbb E(\sigma(X_{s}^{i})^4\varphi_j(X_{s}^{i})^4)ds\\
 & \leqslant &
 N^2TL(m)^2\int_{0}^{T}\mathbb E(\sigma(X_s)^4)ds =
 N^2T^2L(m)^2\int_{-\infty}^{\infty}\sigma(x)^4f_T(x)dx.
\end{eqnarray*}
On the other hand, by Jensen's inequality and Cauchy-Schwarz's inequality,
\begin{eqnarray*}
 \sum_{j = 1}^{m}\mathbb E(A_{j}^{2})
 & \leqslant &
 T\left(\sum_{i\not= k}|R_{i,k}|\right)
 \sum_{j = 1}^{m}\sum_{i\neq k}|R_{i,k}|\int_{0}^{T}\mathbb E(
 \sigma(X_{s}^{i})^2\varphi_j(X_{s}^{i})^2
 \sigma(X_{s}^{k})^2\varphi_j(X_{s}^{k})^2)ds\\
 & \leqslant &
 T\left(\sum_{i\not= k}|R_{i,k}|\right)^2
 \sum_{j = 1}^{m}\int_{0}^{T}\mathbb E(
 \sigma(X_s)^4\varphi_j(X_s)^4)ds\\
 & \leqslant &
 T^2\left(\sum_{i\not= k}|R_{i,k}|\right)^2L(m)^2\int_{-\infty}^{\infty}\sigma(x)^4f_T(x)dx.
\end{eqnarray*}
Therefore,
\begin{displaymath}
\mathbb E(|\widehat{\bf E}_{m}^{*}\widehat{\bf E}_m|^2)
\leqslant
\frac{\mathfrak c_2}{N^2T^2}mL(m)^2\left[1 +
\left(\frac{1}{N}\sum_{i\not= k}|R_{i,k}|\right)^2\right]
\end{displaymath}
with
\begin{displaymath}
\mathfrak c_2 = 2\mathfrak c_1\int_{-\infty}^{\infty}\sigma(x)^4f_T(x)dx.
\end{displaymath}
%

% Subsubsection : Proof of Lemma deviation_probabilities.

%
\subsubsection{Proof of Lemma \ref{deviation_probabilities}}
Let $(\overline\varphi_1,\dots,\overline\varphi_{N_T})$ be the orthonormal family of $\mathbb L^2(I,f_T(x)dx)$ derived from $(\varphi_1,\dots,\varphi_{N_T})$ via Gram-Schmidt's method. Consider also the matrix
\begin{displaymath}
\widehat{\bf G}_m :=
\sum_{i = 1}^{N}\widehat{\bf G}_m(X^i),
\end{displaymath}
where
\begin{displaymath}
\widehat{\bf G}_m(\psi) :=
\frac{1}{NT}
\left(\int_{0}^{T}\overline\varphi_j(\psi(t))\overline\varphi_{\ell}(\psi(t))dt\right)_{j,\ell\in\{1,\dots,m\}}
\textrm{$;$ }
\forall\psi\in\Omega.
\end{displaymath}
The random matrix $\widehat{\bf G}_m(X^i)$ has the same eigenvalues as $N^{-1}\mathbf\Psi_{m}^{-1/2}\widehat{\bf\Psi}_m(X^i)\mathbf\Psi_{m}^{-1/2}$, where
\begin{displaymath}
\widehat{\bf\Psi}_m(\psi) :=
\left(\frac{1}{T}
\int_{0}^{T}\varphi_j(\psi(t))\varphi_{\ell}(\psi(t))dt\right)_{j,\ell\in\{1,\dots,m\}}
\textrm{$;$ }
\forall\psi\in\Omega.
\end{displaymath}
Moreover, for every $\psi\in\Omega$, by Jensen's and Cauchy-Schwarz's inequalities,
\begin{eqnarray}
 \|\widehat{\bf\Psi}_m(\psi)\|_{\rm op}^{2}
 & = &
 \sup_{\|x\|_{2,m} = 1}
 \sum_{\ell = 1}^{m}\left[
 \frac{1}{T}\int_{0}^{T}\left(\sum_{j = 1}^{m}
 \varphi_j(\psi(t))\varphi_{\ell}(\psi(t))x_j\right)dt\right]^2
 \nonumber\\
 & \leqslant &
 \frac{1}{T}\sup_{\|x\|_{2,m} = 1}
 \sum_{\ell = 1}^{m}
 \int_{0}^{T}\varphi_{\ell}(\psi(t))^2\left(\sum_{j = 1}^{m}
 \varphi_j(\psi(t))x_j\right)^2dt
 \nonumber\\
 \label{deviation_probabilities_1}
 & \leqslant &
 \frac{L(m)}{T}\sup_{\|x\|_{2,m} = 1}
 \int_{0}^{T}\left(\sum_{j = 1}^{m}
 \varphi_j(\psi(t))x_j\right)^2dt
 \leqslant
 L(m)^2.
\end{eqnarray}
{\bf Notations:}
\begin{itemize}
 \item The semidefinite order on symmetric matrices is denoted by $\preccurlyeq$.
 \item $\mathbb E_0(.) :=\mathbb E(.)$ and $\mathbb E_i(.) :=\mathbb E(.|\mathcal F_i)$ for every $i\in\{1,\dots,N\}$.
\end{itemize}
First of all, note that
\begin{displaymath}
\|\widehat{\bf G}_m -\mathbf I\|_{\rm op}\leqslant M_N + R_N,
\end{displaymath}
where
\begin{displaymath}
M_N :=
\left\|\sum_{i = 1}^{N}(\widehat{\bf G}_m(X^i) -\mathbb E_{i - 1}(\widehat{\bf G}_m(X^i)))
\right\|_{\rm op}
\quad {\rm and}\quad
R_N :=
\left\|\sum_{i = 1}^{N}(\mathbb E_{i - 1}(
\widehat{\bf G}_m(X^i)) - N^{-1}\mathbf I)\right\|_{\rm op}.
\end{displaymath}
The proof of Lemma \ref{deviation_probabilities} is dissected in four steps. Step 1 deals with a suitable bound on $\mathbb P(M_N >\delta/2)$, $\delta > 0$, step 2 with a suitable bound on $\mathbb P(R_N >\delta/2)$,
\begin{displaymath}
\{\|\mathbf\Psi_{m}^{-1}\|_{\rm op} <\|\widehat{\mathbf\Psi}_{m}^{-1} -
\mathbf\Psi_{m}^{-1}\|_{\rm op}\}\subset\Omega_{m}^{c}
\end{displaymath}
is established in step 3, and the conclusion comes in step 4.
\\
\\
{\bf Step 1.} For any $\delta > 0$, let us establish a suitable bound on $\mathbb P(M_N >\delta)$. For every $i\in\{1,\dots,N\}$, since
\begin{displaymath}
\widetilde{\bf G}_m(X^i) :=
\widehat{\bf G}_m(X^i) -
\mathbb E_{i - 1}(\widehat{\bf G}_m(X^i))
\end{displaymath}
is a symmetric matrix, by Jensen's inequality and by Inequality (\ref{deviation_probabilities_1}),
\begin{eqnarray*}
 (-\widetilde{\bf G}_m(X^i))^2 =\widetilde{\bf G}_m(X^i)^2
 & \preccurlyeq &
 \lambda_{\max}[\widetilde{\bf G}_m(X^i)^2]\mathbf I\\
 & &
 \hspace{1.5cm} =
 \|\widehat{\bf G}_m(X^i) -
 \mathbb E_{i - 1}(\widehat{\bf G}_m(X^i))\|_{\rm op}^{2}\mathbf I\\
 & \preccurlyeq &
 \frac{2}{N^2}(\|\mathbf\Psi_{m}^{-1/2}
 \widehat{\bf\Psi}_m(X^i)\mathbf\Psi_{m}^{-1/2}\|_{\rm op}^{2} +
 \mathbb E_{i - 1}(\|\mathbf\Psi_{m}^{-1/2}
 \widehat{\bf\Psi}_m(X^i)\mathbf\Psi_{m}^{-1/2}\|_{\rm op}^{2}))\mathbf I\\
 & \preccurlyeq &
 \frac{2}{N^2}(\|\widehat{\bf\Psi}_m(X^i)\|_{\rm op}^{2} +
 \mathbb E_{i - 1}(\|\widehat{\bf\Psi}_m(X^i)\|_{\rm op}^{2}))
 \|{\bf\Psi}_{m}^{-1}\|_{\rm op}^{2}\mathbf I
 \preccurlyeq
 \mathbf A_{i}^{2}
\end{eqnarray*}
with
\begin{displaymath}
\mathbf A_{i}^{2} =\frac{4}{N^2}[L(m)(\|{\bf\Psi}_{m}^{-1}\|_{\rm op}\vee 1)]^2\mathbf I.
\end{displaymath}
So, by Azuma's inequality for matrix martingales (see Tropp \cite{TROPP12}, Theorem 7.1),
\begin{displaymath}
\mathbb P\left(\lambda_{\max}\left(
\sum_{i = 1}^{N}\widetilde{\bf G}_m(X^i)\right) >\delta\right)
\leqslant m\exp\left(-\frac{\delta^2}{8\sigma^2}\right)
\end{displaymath}
and
\begin{displaymath}
\mathbb P\left(-\lambda_{\min}\left(
\sum_{i = 1}^{N}\widetilde{\bf G}_m(X^i)\right) >\delta\right) =
\mathbb P\left(\lambda_{\max}\left(
\sum_{i = 1}^{N}(-\widetilde{\bf G}_m(X^i))\right) >\delta\right)
\leqslant m\exp\left(-\frac{\delta^2}{8\sigma^2}\right),
\end{displaymath}
where
\begin{displaymath}
\sigma^2 =
\left\|\sum_{i = 1}^{N}\mathbf A_{i}^{2}\right\|_{\rm op} =
\frac{4}{N}[L(m)(\|{\bf\Psi}_{m}^{-1}\|_{\rm op}\vee 1)]^2.
\end{displaymath}
This leads to
\begin{eqnarray*}
 \mathbb P(M_N >\delta) & = &
 \mathbb P\left(\left\|
 \sum_{i = 1}^{N}\widetilde{\bf G}_m(X^i)\right\|_{\rm op} >\delta\right)\\
 & &
 \hspace{2cm} =
 \mathbb P\left(\max\left\{\lambda_{\max}\left(
 \sum_{i = 1}^{N}\widetilde{\bf G}_m(X^i)\right);
 -\lambda_{\min}\left(
 \sum_{i = 1}^{N}\widetilde{\bf G}_m(X^i)\right)\right\} >\delta\right)\\
 & \leqslant &
 2m\exp\left(-\frac{\delta^2}{8\sigma^2}\right) =
 2m\exp\left[-\frac{\delta^2N}{32[L(m)(\|{\bf\Psi}_{m}^{-1}\|_{\rm op}\vee 1)]^2}\right].
\end{eqnarray*}
{\bf Step 2.} For any $\delta > 0$, let us establish a suitable bound on $\mathbb P(R_N >\delta)$. First of all, let us recall that
\begin{displaymath}
\mathcal I_N =
\{i\in\{1,\dots,N\} : B^i\textrm{ is not independent of }\mathcal F_{i - 1}\}.
\end{displaymath}
For every $i\in\{1,\dots,N\}\backslash\mathcal I_N$, since $(\overline\varphi_1,\dots,\overline\varphi_{N_T})$ is an orthonormal family of $\mathbb L^2(I,f_T(x)dx)$,
\begin{eqnarray*}
 \mathbb E_{i - 1}(\widehat{\bf G}_m(X^i)) & = &
 \mathbb E(\widehat{\bf G}_m(X)) =
 \left(\frac{1}{NT}\int_{0}^{T}\mathbb E(\overline\varphi_j(X_t)
 \overline\varphi_{\ell}(X_t))dt\right)_{j,\ell}\\
 & = &
 \frac{1}{N}
 (\langle\overline\varphi_j,\overline\varphi_{\ell}\rangle_{f_T})_{j,\ell} =
 \frac{1}{N}\mathbf I.
\end{eqnarray*}
Then,
\begin{displaymath}
R_N =
\left\|\sum_{i\in\mathcal I_N}(\mathbb E_{i - 1}(
\widehat{\bf G}_m(X^i)) - N^{-1}\mathbf I)\right\|_{\rm op}.
\end{displaymath}
By Markov's inequality and Jensen's inequality (usual and conditional),
\begin{eqnarray*}
 \mathbb P(R_N >\delta) & \leqslant &
 \frac{\mathbb E(R_{N}^{p})}{\delta^p}
 \leqslant
 \frac{|\mathcal I_N|^{p - 1}}{\delta^p}\sum_{i\in\mathcal I_N}
 \mathbb E(\|\mathbb E_{i - 1}(\widehat{\bf G}_m(X^i) - N^{-1}\mathbf I)\|_{\rm op}^{p})\\
 & \leqslant &
 \frac{|\mathcal I_N|^p}{\delta^p}
 \mathbb E(\|\widehat{\bf G}_m(X) - N^{-1}\mathbf I\|_{\rm op}^{p})\\
 & &
 \hspace{3cm} =
 \frac{|\mathcal I_N|^p}{\delta^pN^p}
 \mathbb E(\|\mathbf\Psi_{m}^{-1/2}\widehat{\bf\Psi}_m(X)\mathbf\Psi_{m}^{-1/2}
 -\mathbf I\|_{\rm op}^{p})\\
 & \leqslant &
 \frac{2^{p - 1}|\mathcal I_N|^p}{\delta^pN^p}[
 \mathbb E(\|\widehat{\bf\Psi}_m(X)\|_{\rm op}^{p})
 \|\mathbf\Psi_{m}^{-1}\|_{\rm op}^{p} + 1]
 \leqslant\frac{2^p|\mathcal I_N|^p}{\delta^pN^p}
 [L(m)(\|\mathbf\Psi_{m}^{-1}\|_{\rm op}\vee 1)]^p.
\end{eqnarray*}
{\bf Step 3.} Now, consider
\begin{displaymath}
\Theta_m :=
\{\|\mathbf\Psi_{m}^{-1}\|_{\rm op} <\|\widehat{\mathbf\Psi}_{m}^{-1} -
\mathbf\Psi_{m}^{-1}\|_{\rm op}\}.
\end{displaymath}
Note that
\begin{displaymath}
\|\widehat{\bf\Psi}_{m}^{-1} -\mathbf\Psi_{m}^{-1}\|_{\rm op} =
\|\mathbf\Psi_{m}^{-\frac{1}{2}}(\widehat{\bf G}_{m}^{-1} -\mathbf I)
\mathbf\Psi_{m}^{-\frac{1}{2}}\|_{\rm op}\leqslant
\|\widehat{\bf G}_{m}^{-1} -\mathbf I\|_{\rm op}\|\mathbf\Psi_{m}^{-1}\|_{\rm op}.
\end{displaymath}
Moreover, as established in Stewart and Sun \cite{SS90}, for every $\mathbf A,\mathbf B\in\mathcal M_d(\mathbb R)$, if $\mathbf A$ is invertible and $\|\mathbf A^{-1}\mathbf B\|_{\rm op} < 1$, then $\mathbf M :=\mathbf A +\mathbf B$ is invertible, and
\begin{displaymath}
\|\mathbf M^{-1} -\mathbf A^{-1}\|_{\rm op}
\leqslant\frac{\|\mathbf B\|_{\rm op}
\|\mathbf A^{-1}\|_{\rm op}^{2}}{1 -\|\mathbf A^{-1}\mathbf B\|_{\rm op}}.
\end{displaymath}
On $\Omega_m$, by applying this result to $\mathbf A =\mathbf I$ and $\mathbf B =\widehat{\bf G}_m -\mathbf I$, $\mathbf A +\mathbf B =\widehat{\bf G}_m$ is invertible and
\begin{displaymath}
\|\widehat{\bf G}_{m}^{-1} -\mathbf I\|_{\rm op}
\leqslant\frac{\|\widehat{\bf G}_m -\mathbf I\|_{\rm op}}{1
-\|\widehat{\bf G}_m -\mathbf I\|_{\rm op}}.
\end{displaymath}
Therefore,
\begin{eqnarray*}
 \Theta_m
 & \subset &
 \{\|\widehat{\bf G}_{m}^{-1} -\mathbf I\|_{\rm op} > 1\}\subset
 \Omega_{m}^{c}\cup
 (\Omega_m\cap\{\|\widehat{\bf G}_{m}^{-1} -\mathbf I\|_{\rm op} > 1\})\\
 & \subset &
 \Omega_{m}^{c}\cup
 \left\{\|\widehat{\bf G}_m -\mathbf I\|_{\rm op}\leqslant\frac{1}{2}
 \textrm{ and }
 \frac{\|\widehat{\bf G}_{m} -\mathbf I\|_{\rm op}}{1
 -\|\widehat{\bf G}_{m} -\mathbf I\|_{\rm op}} > 1\right\}
 =\Omega_{m}^{c}.
\end{eqnarray*}
{\bf Step 4 (conclusion).} For any $\delta > 0$, the two previous steps leads to
\begin{eqnarray}
 \mathbb P(\|\widehat{\bf G}_m -\mathbf I\|_{\rm op} >\delta) & \leqslant &
 \mathbb P\left(
 \left\{M_N >\frac{\delta}{2}\right\}\cup\left\{R_N >\frac{\delta}{2}\right\}\right)
 \nonumber\\
 & \leqslant &
 2m\exp\left[-\frac{\delta^2N}{128[L(m)(\|{\bf\Psi}_{m}^{-1}\|_{\rm op}\vee 1)]^2}\right]
 \nonumber\\
 \label{deviation_probabilities_2}
 & &
 \hspace{4cm} +
 \frac{2^{2p}|\mathcal I_N|^p}{\delta^pN^p}
 [L(m)(\|\mathbf\Psi_{m}^{-1}\|_{\rm op}\vee 1)]^p.
\end{eqnarray}
As established in the beginning of the proof of Comte and Genon-Catalot \cite{CGC20}, Proposition 2.1,
\begin{displaymath}
\Omega_m =
\left\{\sup_{\tau\in\mathcal S_m}
\left|\frac{\|\tau\|_{N}^{2}}{\|\tau\|_{f_T}^{2}} - 1\right|\leqslant\frac{1}{2}
\right\} =
\left\{\|\widehat{\bf G}_m -\mathbf I\|_{\rm op}\leqslant\frac{1}{2}\right\}.
\end{displaymath}
Then, by Inequality (\ref{deviation_probabilities_2}), by Assumptions \ref{condition_m} and \ref{sparsity_condition} (leading to (\ref{sparsity_condition_1})), and since $p\geqslant 12$,
\begin{eqnarray*}
 \mathbb P(\Omega_{m}^{c})
 & \leqslant &
 2m\exp\left(-\frac{\log(NT)}{256\mathfrak c_T(p)T}\right) +
 \frac{2^{3p}(\mathfrak c_T(p)T)^{p/2}}{2^{p/2}}
 \cdot\frac{|\mathcal I_N|^p}{N^{p/2}\log(NT)^{p/2}}\\
 & \leqslant &
 \mathfrak c_1\left(\frac{m}{N^{1 + p/2}} +
 \frac{|\mathcal I_N|^p}{N^{p/2}\log(NT)^{p/2}}\right)
 \leqslant
 \frac{\mathfrak c_1(1 +\mathfrak c_{\ref{sparsity_condition}}(p))}{N^6}
\end{eqnarray*}
where $\mathfrak c_1 > 0$ is a deterministic constant not depending on $m$ and $N$. Moreover, on $\widehat\Lambda_{m}^{c}$ and by Assumption \ref{condition_m},
\begin{displaymath}
[L(m)(\|\mathbf\Psi_{m}^{-1}\|_{\rm op}\vee 1)]^2
\leqslant\frac{\mathfrak c_T(p)}{2}\cdot\frac{NT}{\log(NT)}
\quad {\rm and}\quad
L(m)(\|\widehat{\bf\Psi}_{m}^{-1}\|_{\rm op}\vee 1)
>\mathfrak c_T(p)\frac{NT}{\log(NT)}.
\end{displaymath}
The first inequality implies that
\begin{displaymath}
L(m)\|\mathbf\Psi_{m}^{-1}\|_{\rm op}
\leqslant\frac{\mathfrak c_T(p)}{2}\cdot\frac{NT}{\log(NT)}
\quad {\rm and}\quad
L(m)\leqslant
\frac{\mathfrak c_T(p)}{2}\cdot\frac{NT}{\log(NT)},
\end{displaymath}
and then the second one leads to
\begin{eqnarray*}
 \mathfrak c_T(p)\frac{NT}{\log(NT)} <
 L(m)\|\widehat{\bf\Psi}_{m}^{-1}\|_{\rm op}
 & \leqslant &
 L(m)(
 \|\widehat{\bf\Psi}_{m}^{-1} -\mathbf\Psi_{m}^{-1}\|_{\rm op} +
 \|\mathbf\Psi_{m}^{-1}\|_{\rm op})\\
 & \leqslant &
 L(m)\|\widehat{\bf\Psi}_{m}^{-1} -\mathbf\Psi_{m}^{-1}\|_{\rm op} +
 \frac{\mathfrak c_T(p)}{2}\cdot\frac{NT}{\log(NT)}.
\end{eqnarray*}
Therefore, by step 3,
\begin{eqnarray*}
 \mathbb P(\widehat\Lambda_{m}^{c})
 & \leqslant &
 \mathbb P\left(
 \frac{\mathfrak c_T(p)}{2}\cdot\frac{NT}{\log(NT)}
 \leqslant L(m)
 \|\widehat{\bf\Psi}_{m}^{-1} -\mathbf\Psi_{m}^{-1}\|_{\rm op}\right)\\
 & \leqslant &
 \mathbb P(\|\mathbf\Psi_{m}^{-1}\|_{\rm op} <
 \|\widehat{\bf\Psi}_{m}^{-1} -\mathbf\Psi_{m}^{-1}\|_{\rm op})
 \leqslant\mathbb P(\Omega_{m}^{c})\leqslant
 \frac{\mathfrak c_1(1 +\mathfrak c_{\ref{sparsity_condition}}(p))}{N^6}.
\end{eqnarray*}
%

% Subsection : Proof of Theorem bound_adaptive_estimator.

%
\subsection{Proof of Theorem \ref{bound_adaptive_estimator}}
Let us consider the events
\begin{displaymath}
\Omega_N :=
\bigcap_{m\in\mathcal M_{N}^{+}}\Omega_m
\quad {\rm and}\quad
\Xi_N :=\{\mathcal M_N\subset\widehat{\mathcal M}_N\subset\mathcal M_{N}^{+}\},
\end{displaymath}
where
\begin{displaymath}
\mathcal M_{N}^{+} :=
\left\{m\in\{1,\dots,N_T\} :
[\mathfrak c_{\varphi}^{2}m(\|\mathbf\Psi_{m}^{-1}\|_{\rm op}\vee 1)]^2
\leqslant 4\mathfrak d_T(p)\frac{NT}{\log(NT)}
\right\}.
\end{displaymath}
Moreover, recall that
\begin{displaymath}
\mathcal M_N =
\left\{m\in\{1,\dots,N_T\} :
[\mathfrak c_{\varphi}^{2}m(\|\mathbf\Psi_{m}^{-1}\|_{\rm op}\vee 1)]^2
\leqslant\frac{\mathfrak d_T(p)}{4}\cdot\frac{NT}{\log(NT)}
\right\}\subset\mathcal M_{N}^{+}
\end{displaymath}
and
\begin{displaymath}
\widehat{\mathcal M}_N =
\left\{m\in\{1,\dots,N_T\} :
[\mathfrak c_{\varphi}^{2}m(\|\widehat{\bf\Psi}_{m}^{-1}\|_{\rm op}\vee 1)]^2
\leqslant\mathfrak d_T(p)\frac{NT}{\log(NT)}
\right\}.
\end{displaymath}
The proof of Theorem \ref{bound_adaptive_estimator} relies on the two following lemmas.
%

% Lemma : Bound on the event Xi_N.

%
\begin{lemma}\label{bound_Xi_N}
Under Assumptions \ref{assumption_basis}, \ref{sparsity_condition} and \ref{additional_assumption_basis}, there exists a deterministic constant $\mathfrak c_{\ref{bound_Xi_N}} > 0$, not depending on $m$ and $N$, such that
\begin{displaymath}
\mathbb P(\Xi_{N}^{c})
\leqslant\frac{\mathfrak c_{\ref{bound_Xi_N}}}{N^5}.
\end{displaymath}
\end{lemma}
%

% Lemma : Bound on the empirical process.

%
\begin{lemma}\label{bound_empirical_process}
Consider the empirical process
\begin{displaymath}
\nu_N(\tau) :=
\frac{1}{NT}\sum_{i = 1}^{N}\int_{0}^{T}\sigma(X_{s}^{i})\tau(X_{s}^{i})dB_{s}^{i}
\textrm{ $;$ }
\tau\in\mathcal S_1\cup\dots\cup\mathcal S_{N_T}.
\end{displaymath}
Under Assumptions \ref{assumption_basis} and \ref{additional_assumption_Brownians}, there exist deterministic constants $\kappa_0,\mathfrak c_{\ref{bound_empirical_process}} > 0$, not depending on $N$, such that $\mathfrak c_{\rm cal}\geqslant\kappa_0$ and, for every $m\in\mathcal M_N$,
\begin{displaymath}
\mathbb E\left[\left(
\left[\sup_{\tau\in\mathcal B_{\widehat m,m}}|\nu_N(\tau)|\right]^2
- p(\widehat m,m)\right)_+\mathbf 1_{\Xi_N\cap\Omega_N}\right]
\leqslant
\frac{\mathfrak c_{\ref{bound_empirical_process}}}{NT}
\end{displaymath}
where, for every $m'\in\mathcal M_N$,
\begin{displaymath}
\mathcal B_{m,m'} :=
\{\tau\in\mathcal S_{m\vee m'} :\|\tau\|_{f_T} = 1\}
\quad {\rm and}\quad
p(m,m') :=
\frac{\mathfrak c_{\rm cal}}{8}\cdot\frac{m\vee m'}{NT}
\left(1 +\frac{1}{N}\sum_{i\neq k}|R_{i,k}|\right).
\end{displaymath}
\end{lemma}
%

% Subsubsection : Steps of the proof.

%
\subsubsection{Steps of the proof}
First of all,
\begin{eqnarray}
 \|\widehat b_{\widehat m} - b_I\|_{N}^{2} & = &
 \|\widehat b_{\widehat m} - b_I\|_{N}^{2}\mathbf 1_{\Xi_{N}^{c}} +
 \|\widehat b_{\widehat m} - b_I\|_{N}^{2}\mathbf 1_{\Xi_N}
 \nonumber\\
 \label{bound_adaptive_estimator_1}
 & =: &
 U_1 + U_2.
\end{eqnarray}
Let us find suitable bounds on $\mathbb E(U_1)$ and $\mathbb E(U_2)$.
\begin{itemize}
 \item {\bf Bound on $\mathbb E(U_1)$.} By the definition of $\widehat{\mathcal M}_N$ and by Lemma \ref{bound_variance_remainder},
 \begin{eqnarray*}
  \mathbb E(\|\widehat b_{\widehat m} -\Pi_{N,\widehat m}(b_I)\|_{N}^{2}
  \mathbf 1_{\Xi_{N}^{c}})
  & = &
  \mathbb E[(\widehat{\bf E}_{\widehat m}^{*}\widehat{\bf\Psi}_{\widehat m}^{-1}
  \widehat{\bf E}_{\widehat m})\mathbf 1_{\Xi_{N}^{c}}]\\
  & &
  \hspace{2cm}
  \leqslant
  \mathbb E(\|\widehat{\bf\Psi}_{\widehat m}^{-1}\|_{\rm op}
  |\widehat{\bf E}_{N_T}^{*}\widehat{\bf E}_{N_T}|
  \mathbf 1_{\Xi_{N}^{c}})\\
  & \leqslant &
  \sqrt{
  \mathfrak d_T(p)
  \frac{NT}{\log(NT)}}
  \mathbb E(|\widehat{\bf E}_{N_T}^{*}\widehat{\bf E}_{N_T}|^2)^{1/2}
  \mathbb P(\Xi_{N}^{c})^{1/2}\\
  & \leqslant &
  \frac{\mathfrak c_1N}{\log(NT)}
  \left(1 +\frac{1}{N}\sum_{i\neq k}|R_{i,k}|\right)\mathbb P(\Xi_{N}^{c})^{1/2}
 \end{eqnarray*}
 where $\mathfrak c_1 > 0$ is a deterministic constant not depending on $N$. Then,
 \begin{eqnarray*}
  \mathbb E(U_1) & \leqslant &
  \mathbb E(\|\widehat b_{\widehat m} -\Pi_{N,\widehat m}(b_I)\|_{N}^{2}
  \mathbf 1_{\Xi_{N}^{c}}) +
  \mathbb E(\|b_I\|_{N}^{2}\mathbf 1_{\Xi_{N}^{c}})\\
  & \leqslant &
  \frac{\mathfrak c_1N}{\log(NT)}
  \left(1 +\frac{1}{N}\sum_{i\neq k}|R_{i,k}|\right)\mathbb P(\Xi_{N}^{c})^{1/2} +
  \mathfrak c_2\mathbb P(\Xi_{N}^{c})^{1/2}
 \end{eqnarray*}
 with
 \begin{displaymath}
 \mathfrak c_2 =\left(\int_{-\infty}^{\infty}b_I(x)^4f_T(x)dx\right)^{1/2}.
 \end{displaymath}
 So, by Lemma \ref{deviation_probabilities}, there exists a deterministic constant $\mathfrak c_3 > 0$, not depending on $N$, such that
 \begin{displaymath}
 \mathbb E(U_1)
 \leqslant\frac{\mathfrak c_3}{N}
 \left(1 +\frac{1}{N}\sum_{i\neq k}|R_{i,k}|\right).
 \end{displaymath}
 \item {\bf Bound on $\mathbb E(U_2)$.} Note that
 \begin{eqnarray*}
  U_2 & = &
  \|\widehat b_{\widehat m} - b_I\|_{N}^{2}\mathbf 1_{\Xi_N\cap\Omega_{N}^{c}} +
  \|\widehat b_{\widehat m} - b_I\|_{N}^{2}\mathbf 1_{\Xi_N\cap\Omega_N}\\
  & =: &
  U_{2,1} + U_{2,2}.
 \end{eqnarray*}
 On the one hand, by Lemma \ref{deviation_probabilities}, there exists a deterministic constant $\mathfrak c_4 > 0$, not depending on $N$, such that
 \begin{displaymath}
 \mathbb P(\Xi_N\cap\Omega_{N}^{c})
 \leqslant
 \sum_{m\in\mathcal M_{N}^{+}}\mathbb P(\Omega_{m}^{c})
 \leqslant\frac{\mathfrak c_4}{N^5}.
 \end{displaymath}
 Then, as for $\mathbb E(U_1)$, there exists a deterministic constant $\mathfrak c_5 > 0$, not depending on $N$, such that
 \begin{displaymath}
 \mathbb E(U_{2,1})
 \leqslant\frac{\mathfrak c_5}{N}
 \left(1 +\frac{1}{N}\sum_{i\neq k}|R_{i,k}|\right).
 \end{displaymath}
 On the other hand,
 \begin{displaymath}
 \gamma_N(\tau') -\gamma_N(\tau) =
 \|\tau' - b\|_{N}^{2} -\|\tau - b\|_{N}^{2} - 2\nu_N(\tau' -\tau)
 \end{displaymath}
 for every $\tau,\tau'\in\mathcal S_1\cup\dots\cup\mathcal S_{N_T}$. Moreover, since
 \begin{displaymath}
 \widehat m =
 \arg\min_{m\in\widehat{\mathcal M}_N}\{-\|\widehat b_m\|_{N}^{2} +
 {\rm pen}(m)\} =
 \arg\min_{m\in\widehat{\mathcal M}_N}\{\gamma_N(\widehat b_m) +
 {\rm pen}(m)\},
 \end{displaymath}
 for every $m\in\widehat{\mathcal M}_N$,
 \begin{equation}\label{bound_adaptive_estimator_2}
 \gamma_N(\widehat b_{\widehat m}) + {\rm pen}(\widehat m)
 \leqslant
 \gamma_N(\widehat b_m) + {\rm pen}(m).
 \end{equation}
 On the event $\Xi_N =\{\mathcal M_N\subset\widehat{\mathcal M}_N\subset\mathcal M_{N}^{+}\}$, Inequality (\ref{bound_adaptive_estimator_2}) remains true for every $m\in\mathcal M_N$. Then, on $\Xi_N$, for any $m\in\mathcal M_N$, since $\mathcal S_m +\mathcal S_{\widehat m}\subset\mathcal S_{m\vee\widehat m}$ under Assumption \ref{additional_assumption_basis},
 \begin{eqnarray*}
  \|\widehat b_{\widehat m} - b_I\|_{N}^{2}
  & \leqslant &
  \|\widehat b_m - b_I\|_{N}^{2} +
  2\nu_N(\widehat b_{\widehat m} -\widehat b_m) +
  {\rm pen}(m) - {\rm pen}(\widehat m)\\
  & \leqslant &
  \|\widehat b_m - b_I\|_{N}^{2} +
  \frac{1}{8}\|\widehat b_{\widehat m} -\widehat b_m\|_{f_T}^{2}\\
  & &
  \hspace{1.5cm} +
  8\left(\left[\sup_{\tau\in\mathcal B_{m,\widehat m}}|\nu_N(\tau)|\right]^2
  - p(m,\widehat m)\right)_+ +
  {\rm pen}(m) + 8p(m,\widehat m) - {\rm pen}(\widehat m).
 \end{eqnarray*}
 Since $\|.\|_{f_T}^{2}\mathbf 1_{\Omega_N}\leqslant 2\|.\|_{N}^{2}\mathbf 1_{\Omega_N}$ on $\mathcal S_1\cup\dots\cup\mathcal S_{\max(\mathcal M_{N}^{+})}$, and since $8p(m,\widehat m)\leqslant {\rm pen}(m) + {\rm pen}(\widehat m)$, on $\Xi_N\cap\Omega_N$,
 \begin{displaymath}
 \|\widehat b_{\widehat m} - b_I\|_{N}^{2}
 \leqslant
 3\|\widehat b_m - b_I\|_{N}^{2} + 4{\rm pen}(m) +
 16\left(\left[\sup_{\tau\in\mathcal B_{m,\widehat m}}|\nu_N(\tau)|\right]^2
 - p(m,\widehat m)\right)_+.
 \end{displaymath}
 So, by Lemma \ref{bound_empirical_process},
 \begin{eqnarray*}
  \mathbb E(U_{2,2})
  & \leqslant &
  \min_{m\in\mathcal M_N}
  \{\mathbb E(3\|\widehat b_m - b_I\|_{N}^{2}\mathbf 1_{\Xi_N}) + 4{\rm pen}(m)\} +
  \frac{16\mathfrak c_{\ref{bound_empirical_process}}}{NT}\\
  & \leqslant &
  \mathfrak c_6\min_{m\in\mathcal M_N}
  \left\{\inf_{\tau\in\mathcal S_m}
  \|\tau - b_I\|_{f_T}^{2} +\frac{m}{NT}\left(1 +\frac{1}{N}\sum_{i\neq k}|R_{i,k}|\right)
  \right\} +
  \frac{\mathfrak c_6}{N}
 \end{eqnarray*}
 where $\mathfrak c_6 > 0$ is a deterministic constant not depending on $N$.
\end{itemize}
%

% Subsubsection : Proof of Lemma bound_Xi_N.

%
\subsubsection{Proof of Lemma \ref{bound_Xi_N}}
Note that
\begin{displaymath}
\Xi_{N}^{c} =\{\mathcal M_N\not\subset\widehat{\mathcal M}_N\}\cup
\{\widehat{\mathcal M}_N\not\subset\mathcal M_{N}^{+}\}.
\end{displaymath}
The proof of Lemma \ref{bound_Xi_N} is dissected in three steps. Step 1 deals with a bound on $\mathbb P(\mathcal M_N\not\subset\widehat{\mathcal M}_N)$, step 2 with a bound on $\mathbb P(\|\widehat{\bf\Psi}_m -\mathbf\Psi_m\|_{\rm op} >\delta)$, $\delta > 0$, and step 3 with a bound on $\mathbb P(\widehat{\mathcal M}_N\not\subset\mathcal M_{N}^{+})$.
\\
\\
{\bf Step 1.} On $\{\mathcal M_N\not\subset\widehat{\mathcal M}_N\}$, there exists $m\in\{1,\dots,N_T\}$ such that
\begin{displaymath}
[\mathfrak c_{\varphi}^{2}m(\|\mathbf\Psi_{m}^{-1}\|_{\rm op}\vee 1)]^2
\leqslant\frac{\mathfrak d_T(p)}{4}\cdot\frac{NT}{\log(NT)}
\quad {\rm and}\quad
[\mathfrak c_{\varphi}^{2}m(\|\widehat{\bf\Psi}_{m}^{-1}\|_{\rm op}\vee 1)]^2
>\mathfrak d_T(p)\frac{NT}{\log(NT)}.
\end{displaymath}
The first inequality is equivalent to
\begin{displaymath}
\mathfrak c_{\varphi}^{4}m^2\|\mathbf\Psi_{m}^{-1}\|_{\rm op}^{2}
\leqslant\frac{\mathfrak d_T(p)}{4}\cdot\frac{NT}{\log(NT)}
\quad {\rm and}\quad
\mathfrak c_{\varphi}^{4}m^2
\leqslant\frac{\mathfrak d_T(p)}{4}\cdot\frac{NT}{\log(NT)},
\end{displaymath}
and then the second one leads to
\begin{eqnarray*}
 \mathfrak d_T(p)\frac{NT}{\log(NT)} <
 \mathfrak c_{\varphi}^{4}m^2\|\widehat{\bf\Psi}_{m}^{-1}\|_{\rm op}^{2}
 & \leqslant &
 2\mathfrak c_{\varphi}^{4}m^2(
 \|\widehat{\bf\Psi}_{m}^{-1} -\mathbf\Psi_{m}^{-1}\|_{\rm op}^{2} +
 \|\mathbf\Psi_{m}^{-1}\|_{\rm op}^{2})\\
 & \leqslant &
 2\mathfrak c_{\varphi}^{4}m^2
 \|\widehat{\bf\Psi}_{m}^{-1} -\mathbf\Psi_{m}^{-1}\|_{\rm op}^{2} +
 \frac{\mathfrak d_T(p)}{2}\cdot\frac{NT}{\log(NT)}.
\end{eqnarray*}
So,
\begin{eqnarray*}
 \{\mathcal M_N\not\subset\widehat{\mathcal M}_N\}
 & \subset & \bigcup_{m\in\mathcal M_N}\left\{
 \frac{\mathfrak d_T(p)}{4}\cdot\frac{NT}{\log(NT)}
 \leqslant\mathfrak c_{\varphi}^{4}m^2
 \|\widehat{\bf\Psi}_{m}^{-1} -\mathbf\Psi_{m}^{-1}\|_{\rm op}^{2}\right\}\\
 & \subset & \bigcup_{m\in\mathcal M_N}\{
 \|\mathbf\Psi_{m}^{-1}\|_{\rm op} <
 \|\widehat{\bf\Psi}_{m}^{-1} -\mathbf\Psi_{m}^{-1}\|_{\rm op}\}
 \subset\bigcup_{m\in\mathcal M_N}\Omega_{m}^{c}
\end{eqnarray*}
and, since $\mathfrak d_T(p)/4\leqslant\mathfrak c_T(p)/2$, by Lemma \ref{deviation_probabilities},
\begin{displaymath}
\mathbb P(\mathcal M_N\not\subset\widehat{\mathcal M}_N)
\leqslant
\sum_{m\in\mathcal M_N}\mathbb P(\Omega_{m}^{c})
\leqslant\frac{\mathfrak c_1}{N^5}
\end{displaymath}
where $\mathfrak c_1 > 0$ is a deterministic constant not depending on $N$.
\\
\\
{\bf Step 2.} First of all, note that
\begin{displaymath}
\|\widehat{\bf\Psi}_m -\mathbf\Psi_m\|_{\rm op}\leqslant M_N + R_N,
\end{displaymath}
where
\begin{displaymath}
M_N :=
\frac{1}{N}
\left\|\sum_{i = 1}^{N}(\widehat{\bf\Psi}_m(X^i) -
\mathbb E_{i - 1}(\widehat{\bf\Psi}_m(X^i)))
\right\|_{\rm op}
\textrm{ and }
R_N :=
\frac{1}{N}
\left\|\sum_{i = 1}^{N}(\mathbb E_{i - 1}(
\widehat{\bf\Psi}_m(X^i)) -\mathbb E(\widehat{\bf\Psi}_m(X^i)))\right\|_{\rm op}.
\end{displaymath}
On the one hand, for any $\delta > 0$, let us establish a suitable bound on $\mathbb P(M_N >\delta)$. For every $i\in\{1,\dots,N\}$, since
\begin{displaymath}
\widetilde{\bf\Psi}_m(X^i) :=
\frac{1}{N}(\widehat{\bf\Psi}_m(X^i) -
\mathbb E_{i - 1}(\widehat{\bf\Psi}_m(X^i)))
\end{displaymath}
is a symmetric matrix, by Jensen's inequality and by Inequality (\ref{deviation_probabilities_1}),
\begin{eqnarray*}
 (-\widetilde{\bf\Psi}_m(X^i))^2 =
 \widetilde{\bf\Psi}_m(X^i)^2
 & \preccurlyeq &
 \lambda_{\max}[\widetilde{\bf\Psi}_m(X^i)^2]\mathbf I =
 \frac{1}{N^2}\|\widehat{\bf\Psi}_m(X^i) -
 \mathbb E_{i - 1}(\widehat{\bf\Psi}_m(X^i))\|_{\rm op}^{2}\mathbf I\\
 & \preccurlyeq &
 \frac{2}{N^2}(\|\widehat{\bf\Psi}_m(X^i)\|_{\rm op}^{2} +
 \mathbb E_{i - 1}(\|\widehat{\bf\Psi}_m(X^i)\|_{\rm op}^{2}))\mathbf I
 \preccurlyeq
 \mathbf A_{i}^{2}
\end{eqnarray*}
with
\begin{displaymath}
\mathbf A_{i}^{2} :=\frac{4L(m)^2}{N^2}\mathbf I.
\end{displaymath}
So, by Azuma's inequality for matrix martingales (see Tropp \cite{TROPP12}, Theorem 7.1),
\begin{eqnarray*}
 \mathbb P(M_N >\delta) & = &
 \mathbb P\left(\left\|
 \sum_{i = 1}^{N}\widetilde{\bf\Psi}_m(X^i)\right\|_{\rm op} >\delta\right)\\
 & &
 \hspace{2cm} =
 \mathbb P\left(\max\left\{\lambda_{\max}\left(
 \sum_{i = 1}^{N}\widetilde{\bf\Psi}_m(X^i)\right) ;
 -\lambda_{\min}\left(
 \sum_{i = 1}^{N}\widetilde{\bf\Psi}_m(X^i)\right)\right\} >\delta\right)\\
 & \leqslant &
 2m\exp\left(-\frac{\delta^2N}{32L(m)^2}\right).
\end{eqnarray*}
On the other hand, let us establish a suitable bound on $\mathbb P(R_N >\delta)$. By the definition of $\mathcal I_N$,
\begin{displaymath}
R_N =\frac{1}{N}
\left\|\sum_{i\in\mathcal I_N}(\mathbb E_{i - 1}(
\widehat{\bf\Psi}_m(X^i)) -\mathbb E(\widehat{\bf\Psi}_m(X^i)))\right\|_{\rm op}.
\end{displaymath}
Then, by Markov's inequality and Jensen's inequality (usual and conditional),
\begin{eqnarray*}
 \mathbb P(R_N >\delta) & \leqslant &
 \frac{\mathbb E(R_{N}^{p})}{\delta^p}\leqslant
 \frac{|\mathcal I_N|^{p - 1}}{\delta^pN^p}
 \sum_{i\in\mathcal I_N}\mathbb E[\|
 \mathbb E_{i - 1}[\widehat{\bf\Psi}_m(X^i) -
 \mathbb E(\widehat{\bf\Psi}_m(X^i))]\|_{\rm op}^{p}]\\
 & \leqslant &
 \frac{|\mathcal I_N|^p}{\delta^pN^p}
 \mathbb E(\|\widehat{\bf\Psi}_m(X) -
 \mathbb E(\widehat{\bf\Psi}_m(X))\|_{\rm op}^{p})
 \leqslant\frac{2^p|\mathcal I_N|^p}{\delta^pN^p}L(m)^p.
\end{eqnarray*}
Therefore,
\begin{eqnarray*}
 \mathbb P(\|\widehat{\bf\Psi}_m -\mathbf\Psi_m\|_{\rm op} >\delta) & \leqslant &
 \mathbb P\left(
 \left\{M_N >\frac{\delta}{2}\right\}\cup\left\{R_N >\frac{\delta}{2}\right\}\right)\\
 & \leqslant &
 2m\exp\left(-\frac{\delta^2N}{128L(m)^2}\right) +
 \frac{2^{2p}|\mathcal I_N|^p}{\delta^pN^p}L(m)^p.
\end{eqnarray*}
{\bf Step 3.} On $\{\widehat{\mathcal M}_N\not\subset\mathcal M_{N}^{+}\}$, there exists $m\in\{1,\dots,N_T\}$ such that
\begin{displaymath}
[\mathfrak c_{\varphi}^{2}m(\|\widehat{\bf\Psi}_{m}^{-1}\|_{\rm op}\vee 1)]^2
\leqslant\mathfrak d_T(p)\frac{NT}{\log(NT)}
\quad {\rm and}\quad
[\mathfrak c_{\varphi}^{2}m(\|\mathbf\Psi_{m}^{-1}\|_{\rm op}\vee 1)]^2
> 4\mathfrak d_T(p)\frac{NT}{\log(NT)}.
\end{displaymath}
The first inequality is equivalent to
\begin{displaymath}
\mathfrak c_{\varphi}^{4}m^2\|\widehat{\bf\Psi}_{m}^{-1}\|_{\rm op}^{2}
\leqslant\mathfrak d_T(p)\frac{NT}{\log(NT)}
\quad {\rm and}\quad
\mathfrak c_{\varphi}^{4}m^2
\leqslant\mathfrak d_T(p)\frac{NT}{\log(NT)},
\end{displaymath}
and then the second one leads to
\begin{eqnarray*}
 4\mathfrak d_T(p)\frac{NT}{\log(NT)} <
 \mathfrak c_{\varphi}^{4}m^2\|\mathbf\Psi_{m}^{-1}\|_{\rm op}^{2}
 & \leqslant &
 2\mathfrak c_{\varphi}^{4}m^2(
 \|\mathbf\Psi_{m}^{-1} -\widehat{\bf\Psi}_{m}^{-1}\|_{\rm op}^{2} +
 \|\widehat{\bf\Psi}_{m}^{-1}\|_{\rm op}^{2})\\
 & \leqslant &
 2\mathfrak c_{\varphi}^{4}m^2
 \|\mathbf\Psi_{m}^{-1} -\widehat{\bf\Psi}_{m}^{-1}\|_{\rm op}^{2} +
 2\mathfrak d_T(p)\frac{NT}{\log(NT)}.
\end{eqnarray*}
Moreover, for every $m\in\{1,\dots,N_T\}$,
\begin{eqnarray*}
 & &
 \{\|\mathbf\Psi_{m}^{-1} -\widehat{\bf\Psi}_{m}^{-1}\|_{\rm op} >
 \|\widehat{\bf\Psi}_{m}^{-1}\|_{\rm op}\}\\
 & &
 \hspace{3cm}\subset
 \left\{\|\widehat{\bf\Psi}_{m}^{-1/2}\mathbf\Psi_m\widehat{\bf\Psi}_{m}^{-1/2} -
 \mathbf I\|_{\rm op} >\frac{1}{2}\right\}\subset
 \left\{\|\widehat{\bf\Psi}_m -\mathbf\Psi_m\|_{\rm op} >\frac{1}{2}
 \|\widehat{\bf\Psi}_{m}^{-1}\|_{\rm op}^{-1}\right\}
\end{eqnarray*}
by interchanging $\widehat{\bf\Psi}_m$ and $\mathbf\Psi_m$ in the proof of Comte and Genon-Catalot \cite{CGC19}, Proposition 4.(ii). So,
\begin{eqnarray*}
 \{\widehat{\mathcal M}_N\not\subset\mathcal M_{N}^{+}\}
 & \subset &
 \bigcup_{\mathfrak c_{\varphi}^{4}m^2
 \leqslant\mathfrak d_T(p)NT/\log(NT)}
 \left\{2\mathfrak c_{\varphi}^{4}m^2
 \|\widehat{\bf\Psi}_{m}^{-1}\|_{\rm op}^{2}
 < 2\mathfrak d_T(p)\frac{NT}{\log(NT)}
 \leqslant 2\mathfrak c_{\varphi}^{4}m^2
 \|\mathbf\Psi_{m}^{-1} -\widehat{\bf\Psi}_{m}^{-1}\|_{\rm op}^{2}\right\}\\
 & \subset &
 \bigcup_{\mathfrak c_{\varphi}^{4}m^2
 \leqslant\mathfrak d_T(p)NT/\log(NT)}\left\{
 \|\widehat{\bf\Psi}_m -\mathbf\Psi_m\|_{\rm op} >
 \frac{m}{2}\sqrt{
 \frac{\log(NT)}{\mathfrak d_T(p)NT}}\right\}
\end{eqnarray*}
and, by the previous step, Assumptions \ref{sparsity_condition} and \ref{additional_assumption_basis}, and since $p\geqslant 12$,
\begin{eqnarray*}
 \mathbb P(\widehat{\mathcal M}_N\not\subset\mathcal M_{N}^{+}) & \leqslant &
 \sum_{\mathfrak c_{\varphi}^{4}m^2\leqslant\mathfrak d_T(p)NT/\log(NT)}
 \left(2m\exp\left(-\frac{N}{512L(m)^2}\cdot\frac{m^2\log(NT)}{\mathfrak d_T(p)NT}\right)
 \right.\\
 & &
 \hspace{5cm}\left. +
 \frac{2^{3p}|\mathcal I_N|^p}{N^p}\cdot
 \frac{\mathfrak d_T(p)^{p/2}(NT)^{p/2}}{m^p\log(NT)^{p/2}}L(m)^p\right)\\
 & \leqslant &
 \sum_{\mathfrak c_{\varphi}^{4}m^2\leqslant\mathfrak d_T(p)NT/\log(NT)}
 \left(2m\exp\left(-\frac{1}{512\mathfrak c_{\varphi}^{4}}\cdot\frac{\log(NT)}{\mathfrak d_T(p)T}\right)
 \right.\\
 & &
 \hspace{5cm}\left. +
 \frac{2^{3p}|\mathcal I_N|^p}{N^{p/2}}\cdot
 \frac{\mathfrak d_T(p)^{p/2}T^{p/2}}{\log(NT)^{p/2}}\mathfrak c_{\varphi}^{2p}\right)\\
 & \leqslant &
 \sum_{\mathfrak c_{\varphi}^{4}m^2\leqslant\mathfrak d_T(p)NT/\log(NT)}
 \left(\frac{2m}{N^{1 + p/2}} +
 \frac{2^{3p}\mathfrak c_{\varphi}^{2p}
 \mathfrak c_{\ref{sparsity_condition}}(p)}{N^6}
 \mathfrak d_T(p)^{p/2}T^{p/2}\right)
 \leqslant\frac{\mathfrak c_2}{N^5}
\end{eqnarray*}
where $\mathfrak c_2 > 0$ is a deterministic constant not depending on $N$.
%

% Subsubsection : Proof of Lemma bound_empirical_process.

%
\subsubsection{Proof of Lemma \ref{bound_empirical_process}}
The proof of Lemma \ref{bound_empirical_process} is dissected in two steps.
\\
\\
{\bf Step 1.} Consider $\tau\in\mathcal S_1\cup\dots\cup\mathcal S_{N_T}$ and the martingale $(M_N(\tau)_t)_{t\in [0,T]}$ defined by
\begin{displaymath}
M_N(\tau)_t :=
\sum_{i = 1}^{N}\int_{0}^{t}\sigma(X_{s}^{i})\tau(X_{s}^{i})dB_{s}^{i}
\textrm{ $;$ }
\forall t\in [0,T].
\end{displaymath}
Note that $\nu_N(\tau) = M_N(\tau)_T/(NT)$. Since $d\langle B^i,B^k\rangle_t = R_{i,k}dt$ for every $i,k\in\{1,\dots,N\}$,
\begin{eqnarray*}
 \langle M_N(\tau)\rangle_T & = &
 \sum_{i,k = 1}^{N}R_{i,k}\int_{0}^{T}
 \sigma(X_{t}^{i})\sigma(X_{t}^{k})\tau(X_{t}^{i})\tau(X_{t}^{k})dt\\
 & &
 \hspace{3cm} =
 \int _{0}^{T}(\sigma(X_{t}^{i})\tau(X_{t}^{i}))_{i}^{*}\times R\times
 (\sigma(X_{t}^{i})\tau(X_{t}^{i}))_idt\\
 & \leqslant &
 \|R\|_{\rm op}\int_{0}^{T}
 \|(\sigma(X_{t}^{i})\tau(X_{t}^{i}))_{1\leqslant i\leqslant N}\|_{2,N}^{2}dt\\
 &\leqslant &
 \|R\|_{\rm op}\|\sigma\|_{\infty}^{2}
 \int_{0}^{T}\left(\sum_{i = 1}^{N}\tau(X_{t}^{i})^2\right)dt =
 NT\|R\|_{\rm op}\|\sigma\|_{\infty}^{2}\|\tau\|_{N}^{2}.
\end{eqnarray*}
Then, by Assumption \ref{additional_assumption_Brownians} and Bernstein's inequality for local martingales (see Revuz and Yor \cite{RY99}, p. 153), for any $\varepsilon,\upsilon > 0$,
\begin{eqnarray*}
 \mathbb P(\nu_N(\tau)\geqslant\varepsilon,\|\tau\|_{N}^{2}\leqslant\upsilon^2)
 & \leqslant &
 \mathbb P(M_N(\tau)_{T}^{*}\geqslant NT\varepsilon,
 \langle M_N(\tau)\rangle_T\leqslant
 NT\upsilon^2\|R\|_{\rm op}\|\sigma\|_{\infty}^{2})\\
 & \leqslant &
 \exp\left(
 -\frac{NT\varepsilon^2}{2\upsilon^2\|\sigma\|_{\infty}^{2}\mathfrak r}\right)
 \quad {\rm with}\quad
 \mathfrak r = 1 +
 \mathfrak m_{\ref{additional_assumption_Brownians}}.
\end{eqnarray*}
Since this bound remains true by replacing $\tau$ by $-\tau$,
\begin{eqnarray*}
 \mathbb P(|\nu_N(\tau)|\geqslant\varepsilon,\|\tau\|_{N}^{2}\leqslant\upsilon^2)
 & = &
 \mathbb P(\nu_N(\tau)\geqslant\varepsilon,\|\tau\|_{N}^{2}\leqslant\upsilon^2) +
 \mathbb P(\nu_N(-\tau)\geqslant\varepsilon,\|\tau\|_{N}^{2}\leqslant\upsilon^2)\\
 & \leqslant &
 2\exp\left(
 -\frac{NT\varepsilon^2}{2\upsilon^2\|\sigma\|_{\infty}^{2}\mathfrak r}\right).
\end{eqnarray*}
{\bf Step 2.} By using the bound of step 1 and by following the pattern of the proof of Baraud {\it et al.} \cite{BCV01}, Proposition 6.1, the purpose of this step is to find a suitable bound on
\begin{displaymath}
\mathbb E\left[\left(
\left[\sup_{\tau\in\mathcal B_{m,m'}}|\nu_N(\tau)|\right]^2
- p(m,m')\right)_+\mathbf 1_{\Xi_N\cap\Omega_N}\right]
\textrm{ $;$ }
m,m'\in\mathcal M_N.
\end{displaymath}
Consider $\delta_0\in (0,1)$ and let $(\delta_n)_{n\in\mathbb N^*}$ be the real sequence defined by
\begin{displaymath}
\delta_n :=\delta_02^{-n}
\textrm{ $;$ }
\forall n\in\mathbb N^*.
\end{displaymath}
Since $\mathcal S_{m\vee m'}$ is a vector subspace of $\mathbb L^2(I,f_T(x)dx)$ of dimension $m\vee m'$, by Lorentz {\it et al.} \cite{LGM96}, Chapter 15, Proposition 1.3, for any $n\in\mathbb N$, there exists $T_n\subset\mathcal B_{m,m'}$ such that $|T_n|\leqslant (3/\delta_n)^{m\vee m'}$ and, for any $\tau\in\mathcal B_{m,m'}$,
\begin{displaymath}
\exists f_n\in T_n :
\|\tau - f_n\|_{f_T}\leqslant\delta_n.
\end{displaymath}
In particular, note that
\begin{displaymath}
\tau =
f_0 +\sum_{n = 1}^{\infty}(f_n - f_{n - 1}).
\end{displaymath}
Then, for any sequence $(\Delta_n)_{n\in\mathbb N}$ of elements of $(0,\infty)$ such that $\Delta =\sum_{n\in\mathbb N}\Delta_n <\infty$,
\begin{eqnarray*}
 & & \left\{\left[
 \sup_{\tau\in\mathcal B_{m,m'}}|\nu_N(\tau)|\right]^2 >
 \Delta^2\right\}\\
 & &
 \hspace{2cm} =
 \left\{\exists (f_n)_{n\in\mathbb N}\in\prod_{n = 0}^{\infty}T_n :
 |\nu_N(f_0)| +\sum_{n = 1}^{\infty}
 |\nu_N(f_n - f_{n - 1})| >\Delta\right\}\\
 & &
 \hspace{2cm}\subset
 \left\{\exists (f_n)_{n\in\mathbb N}\in\prod_{n = 0}^{\infty}T_n :
 |\nu_N(f_0)| >\Delta_0\textrm{ or }[\exists n\in\mathbb N^* :
 |\nu_N(f_n - f_{n - 1})| >\Delta_n]\right\}\\
 & &
 \hspace{2cm}\subset
 \bigcup_{f_0\in T_0}\{|\nu_N(f_0)| >\Delta_0\}\cup
 \bigcup_{n = 1}^{\infty}\bigcup_{(f_{n - 1},f_n)\in\mathbb T_n}
 \{|\nu_N(f_n - f_{n - 1})| >\Delta_n\}
\end{eqnarray*}
with $\mathbb T_n = T_{n - 1}\times T_n$ for every $n\in\mathbb N^*$. Moreover, $\|f_0\|_{f_T}^{2}\leqslant 1$,
\begin{displaymath}
\|f_n - f_{n - 1}\|_{f_T}^{2}
\leqslant
2\delta_{n - 1}^{2} + 2\delta_{n}^{2} =\frac{5}{2}\delta_{n - 1}^{2}
\textrm{ $;$ }
\forall n\in\mathbb N^*,
\end{displaymath}
and $\|.\|_{N}^{2}\mathbf 1_{\Omega_N}\leqslant 3/2\|.\|_{f_T}^{2}\mathbf 1_{\Omega_N}$ on $\mathcal S_1\cup\dots\cup\mathcal S_{\max(\mathcal M_{N}^{+})}$. So, by step 1,
\begin{eqnarray}
 & &
 \mathbb P\left(\left\{\left[
 \sup_{\tau\in\mathcal B_{m,m'}}|\nu_N(\tau)|\right]^2 >\Delta^2\right\}
 \cap\Xi_N\cap\Omega_N\right)
 \nonumber\\
 & &
 \hspace{3cm}\leqslant
 2\sum_{f_0\in T_0}\exp\left(
 -\frac{NT\Delta_{0}^{2}}{3\|f_0\|_{f_T}^{2}\|\sigma\|_{\infty}^{2}\mathfrak r}\right)
 \nonumber\\
 & &
 \hspace{6cm}
 + 2\sum_{n = 1}^{\infty}\sum_{(f_{n - 1},f_n)\in\mathbb T_n}
 \exp\left(
 -\frac{NT\Delta_{n}^{2}}{3\|f_n - f_{n - 1}\|_{f_T}^{2}\|\sigma\|_{\infty}^{2}\mathfrak r}\right)
 \nonumber\\
 \label{bound_empirical_process_1}
 & &
 \hspace{3cm}\leqslant
 2\exp\left(h_0 -\frac{NT\Delta_{0}^{2}}{3\|\sigma\|_{\infty}^{2}
 \mathfrak r}\right) +
 2\sum_{n = 1}^{\infty}\exp\left(h_{n - 1} + h_n -
 \frac{NT\Delta_{n}^{2}}{15/2\delta_{n - 1}^{2}\|\sigma\|_{\infty}^{2}\mathfrak r}\right)
\end{eqnarray}
with $h_n =\log(|T_n|)$ for every $n\in\mathbb N$. Now, let us take $\Delta_0$ such that
\begin{displaymath}
h_0 -\frac{NT\Delta_{0}^{2}}{3\|\sigma\|_{\infty}^{2}\mathfrak r} =
-(m\vee m' + x)
\quad {\rm with}\quad x > 0,
\end{displaymath}
which leads to
\begin{displaymath}
\Delta_0 =
\left[\frac{3\|\sigma\|_{\infty}^{2}
\mathfrak r}{NT}(m\vee m' + x + h_0)\right]^{1/2}
\leqslant
\left[\frac{3\|\sigma\|_{\infty}^{2}
\mathfrak r}{NT}(1 + h_0)(m\vee m' + x)\right]^{1/2},
\end{displaymath}
and for every $n\in\mathbb N^*$, let us take $\Delta_n$ such that
\begin{displaymath}
h_{n - 1} + h_n -
\frac{NT\Delta_{n}^{2}}{15/2\delta_{n - 1}^{2}\|\sigma\|_{\infty}^{2}
\mathfrak r} = -(m\vee m' + x + n),
\end{displaymath}
which leads to
\begin{eqnarray*}
 \Delta_n & = &
 \left[\frac{15/2\delta_{n - 1}^{2}\|\sigma\|_{\infty}^{2}\mathfrak r}{NT}
 (m\vee m' + x + n + h_{n - 1} + h_n)\right]^{1/2}\\
 & \leqslant &
 \left[\frac{15/2\delta_{n - 1}^{2}\|\sigma\|_{\infty}^{2}\mathfrak r}{NT}
 (1 + n +\lambda_n)(m\vee m' + x)\right]^{1/2}
 \quad {\rm with}\quad
 \lambda_n = 2\log\left(\frac{3}{\delta_0}\right) + (2n - 1)\log(2).
\end{eqnarray*}
For this appropriate sequence $(\Delta_n)_{n\in\mathbb N}$,
\begin{displaymath}
\mathbb P\left(\left\{\left[
\sup_{\tau\in\mathcal B_{m,m'}}|\nu_N(\tau)|\right]^2 >\Delta^2\right\}
\cap\Xi_N\cap\Omega_N\right)\leqslant
2e^{-x}e^{-(m\vee m')}\left(1 +\sum_{n = 1}^{\infty}e^{-n}\right)\leqslant
3.2e^{-x}e^{-(m\vee m')}
\end{displaymath}
by Inequality (\ref{bound_empirical_process_1}), and
\begin{eqnarray*}
 \Delta^2
 & \leqslant &
 \frac{3\|\sigma\|_{\infty}^{2}\mathfrak r}{NT}
 \left[(1 + h_0)^{1/2}(m\vee m' + x)^{1/2} +
 \sqrt{\frac{5}{2}}\sum_{n = 1}^{\infty}\delta_{n - 1}
 (1 + n +\lambda_n)^{1/2}(m\vee m' + x)^{1/2}\right]^2\\
 & \leqslant &
 \frac{3\|\sigma\|_{\infty}^{2}\mathfrak r}{NT}
 (m\vee m' + x)\delta^2
\end{eqnarray*}
with
\begin{displaymath}
\delta =
(1 + h_0)^{1/2} +
\sqrt{\frac{5}{2}}\sum_{n = 1}^{\infty}\delta_{n - 1}
(1 + n +\lambda_n)^{1/2} <\infty.
\end{displaymath}
Then,
\begin{displaymath}
\mathbb P\left(\left[
\sup_{\tau\in\mathcal B_{m,m'}}|\nu_N(\tau)|\right]^2 -
\frac{8\kappa_0}{\mathfrak c_{\rm cal}TR_N}p(m,m') >
\frac{\kappa_0}{NT}x\right)
\leqslant 3.2e^{-x}e^{-(m\vee m')}
\end{displaymath}
with
\begin{displaymath}
\kappa_0 = 3\|\sigma\|_{\infty}^{2}\mathfrak r\delta^2
\quad {\rm and}\quad
R_N = 1 +\frac{1}{N}\sum_{i\neq k}|R_{i,k}|.
\end{displaymath}
So, by taking $\mathfrak c_{\rm cal}\geqslant 8\kappa_0T^{-1} > 8\kappa_0(TR_N)^{-1}$ and $y =\kappa_0x(NT)^{-1}$,
\begin{displaymath}
\mathbb P\left(\left[
\sup_{\tau\in\mathcal B_{m,m'}}|\nu_N(\tau)|\right]^2 - p(m,m') > y\right)
\leqslant 3.2e^{-NTy/\kappa_0}e^{-(m\vee m')}.
\end{displaymath}
Therefore,
\begin{eqnarray*}
 \mathbb E\left[\left(\left[
 \sup_{\tau\in\mathcal B_{m,m'}}|\nu_N(\tau)|\right]^2 - p(m,m')\right)_+\right]
 & = &
 \int_{0}^{\infty}\mathbb P\left(\left[
 \sup_{\tau\in\mathcal B_{m,m'}}|\nu_N(\tau)|\right]^2 - p(m,m') > y\right)dy\\
 & \leqslant &
 3.2\kappa_0\frac{e^{-(m\vee m')}}{NT}.
\end{eqnarray*}
A union-bound allows to conclude.
%

% References.

%

%
\end{document}